\documentclass{article}
\usepackage{amssymb}
\usepackage{amsmath}
\newtheorem{Theorem}{Theorem}[section]

\newtheorem{Proposition}{Proposition}[section]
\newtheorem{Lemma}{Lemma}[section]
\newtheorem{Corollary}{Corollary}[section]

\newtheorem{Remark}{Remark}[section]
\newtheorem{Example}{Example}[section]

\newcommand{\rec}[1]{{(\ref{#1})}}
\def \R{I\!\!R}

\def\11{1\!\!1}

\newcommand{\rsen}   {{risk-sensitive }}

\newcommand{\ftpm} {{${\cal F}_t$-progressively measurable }}
\def\rec#1{{(\ref{#1})}}
\newcommand{\expect} {{\rm E}}

\newcommand{\jjep} {J^\eps}
\newcommand{\er}[1]{\hbox{(\ref{#1})}}
 \newcommand{\vvep} {V^\eps}
\newcommand{\eps}       {\varepsilon}

\newcommand{\ba}{\begin{array}}
\newcommand{\ea}{\end{array}}

\newcommand{\eptwogam} {\frac{\eps}{2\gamma^2}}

\newcommand{\invtwogam} {\frac{1}{ 2\gamma^2}}

\newcommand{\gamovertwo} {\frac{\gamma^2}{ 2}}

\begin{document}

\title{Uniqueness  Results for Second
Order Bellman-Isaacs Equations under
 Quadratic Growth Assumptions and Applications}
\author{Francesca Da Lio\thanks{Department of Mathematics, University of Padova,
    via Belzoni, 7 35131 Padova, Italy (dalio@math.unipd.it).}\and
Olivier Ley\thanks {Laboratoire de Math\'ematiques et Physique
Th\'eorique. UMR CNRS 6083. Universit\'e de Tours.  Facult\'e des
Sciences et Techniques, Parc de Grandmont, 37200 Tours, France (ley@lmpt.univ-tours.fr).}}

\date{}

\maketitle

\begin{abstract}
In this paper, we prove a comparison result between semicontinuous viscosity sub and supersolutions
growing at most quadratically   of  second-order degenerate parabolic Hamilton-Jacobi-Bellman and
Isaacs equations. As an application, we characterize the value function of a finite horizon
stochastic control problem with unbounded controls as the unique viscosity solution of the
corresponding dynamic programming equation.
\end{abstract}


\pagestyle{myheadings} \thispagestyle{plain} \markboth{F. DA LIO
AND O.LEY}{UNIQUENESS RESULTS FOR BELLMAN-ISAACS EQUATIONS }
\section{Introduction}

In this paper, we are interested in the second-order equation
\begin{equation}\label{hjbintr}
\left\{\begin{array}{ll}
\displaystyle\frac{\partial w}{\partial t}
+ H(x,t,Dw, D^2 w)
+ G(x,t,Dw, D^2 w) = 0 & \mbox{in  $\R^N\times (0,T),$}\\[3mm]
w (x,0)=\psi (x)& \mbox{in  $\R^N,$}\end{array}\right.\
\end{equation}
where $N\geq 1,$ $T>0,$ the unknown $w$ is a real-valued function defined in $\R^N\times[0,T],$
$Dw$ and $D^2w$ denote respectively its gradient and Hessian matrix
and $\psi$ is a given initial condition.
The Hamiltonians $H,G~\colon\R^N\times[0,T]\times \R^N\times{\cal{S}}_{N}(\R)\to\R$
are continuous in all their variables and have the form
\begin{eqnarray} \label{expression-H-intr}
H(x,t,p,X) = \mathop{\rm inf}_{\alpha\in A}
\left\{ \langle b(x,t,\alpha),  p\rangle + \ell(x,t,\alpha)
-{\rm Tr}\left[ \sigma(x,t, \alpha)\sigma^T(x,t,\alpha)
 X \right]
\right\}
\end{eqnarray}
and
\begin{eqnarray} \label{expression-G-intr}
G(x,t, p,X) =
\mathop{\rm sup}_{\beta\in B}
\left\{ -\langle g(x,t, \beta ),  p \rangle
- f (x,t, \beta) -{\rm Tr}\left[ c(x,t, \beta)c^T(x,t,\beta)X \right] \right\}.
\end{eqnarray}
This kind of equation is of particular
interest for applications since it relies on differential game theory (Isaacs equations)
or on deterministic and stochastic control problems when either
$H\equiv 0$ or $G\equiv 0$ (Hamilton-Jacobi-Bellman equations).

Notations and precise assumptions are given in Section \ref{sec-1-comp} but
we point out that we allow one of the control set $A$ or $B$ to be unbounded
and the solutions to (\ref{hjbintr}) may have quadratic growth.
Our model case is the well-known stochastic linear
quadratic problem. We refer to Bensoussan \cite{bensoussan82},
Fleming and Rishel \cite{fr75}, Fleming and Soner \cite{fs93}, {\O}ksendal \cite{oksendal98},
Yong and Zhou \cite{yz99} and the references therein for an overview and to
Examples \ref{LQ-deterministe} and \ref{stochastic-LQ} below. This problem can be described as follows.
Let $(\Omega, {\cal{F}}, ({\cal{F}}_t)_{t\geq 0}, {P})$ be a filtered probability
space and let $(W_t)_{t}$ be a ${\cal{F}}_t$-adapted standard $M$-Brownian motion.
The control set is $A=\R^k$ for some $k>0$ and we consider the linear stochastic differential equation
\begin{eqnarray*}
\left\{\begin{array}{l}
\displaystyle dX_s=[A(s)X_s+B(s)\alpha_s]ds+[C(s)X_s+D(s)]dW_s, \ \ \ {\rm for} \ t\leq s \leq T\\
 X_t=x
\end{array}\right.
\end{eqnarray*}
where $\alpha_s \in{\cal{A}}_t,$ the set of $A$-valued  ${\cal{F}}_t$-progressively measurable
controls and the adapted process $X_s$ is the solution.
The linear quadratic problem consists in minimizing the quadratic cost
\begin{equation}\label{fctval-LQintr}
 V(x,t)=\inf_{\alpha_s \in {\cal A}_t}  \mbox{E}\{\int_t^T
[\langle X_s, Q(s)X_s\rangle+ R|\alpha_s|^2 ]\,ds+\langle X_T, SX_T\rangle \},
\end{equation}
where $A(\cdot),B(\cdot),C(\cdot), D(\cdot), Q(\cdot)$ and $S$ are deterministic
matrix-valued functions of suitable size and $R>0$ to simplify the presentation.
The Hamilton-Jacobi equation associated to this problem is
\begin{equation}\label{LQ-eq-intro}
\left\{\begin{array}{l}
\displaystyle -\frac{\partial w}{\partial t}
-\langle A(t)x,Dw \rangle -\langle x,Q(t)x \rangle
+ \frac{1}{4R}|B(t)^T Dw|^2
-{\rm Tr}\left[ a(x,t) D^2w \right] =0 \\[2mm]
w (x,T)= \langle x,Sx \rangle,
\end{array}\right.
\end{equation}
where $a(x,t)= (C(t)x+D(t))(C(t)x+D(t))^T/2.$
Note that this equation is of type (\ref{hjbintr}) (with $G\equiv 0$) since
\begin{equation}
\frac{1}{4R}|B(t)^T Dw|^2 =
\mathop{\rm sup}_{\alpha\in \R^k}
\Big\{ -\langle B(t) \alpha ,  Dw \rangle
-R |\alpha|^2 \Big\}.
\end{equation}
In this paper, we are concerned with two issues about this problem.

The first question relies on the partial differential equation (\ref{LQ-eq-intro}).
We note that the quadratic cost with unbounded controls leads
to a quadratic term with respect to the gradient variable. From
the terminal condition, we expect the solutions have quadratic growth. Moreover,
the diffusion matrix may be degenerate. Therefore, we cannot hope
to obtain smooth solutions in general. We need to consider weak solutions,
namely viscosity solutions. (We refer the reader who is not
familiar with this notion of solutions to Crandall, Ishii and Lions \cite{cil92},
Fleming and Soner \cite{fs93}, Bardi and Capuzzo Dolcetta \cite{bcd97} and Barles \cite{barles94}
and all the references therein). We obtain the existence of a unique continuous
viscosity solution for (\ref{LQ-eq-intro}) and for a large class of equations
of type (\ref{hjbintr}) (see Theorem \ref{thm-unicite} and Corollary \ref{thm-existence}).

We point out that the results obtained in this paper are beyond the classical comparison
results for viscosity solutions (see e.g. \cite{cil92}) because of the growth both of the
solutions and the Hamiltonians. In fact, most of the comparison results in the literature
require that either the solutions are uniformly continuous or the
Hamiltonian is uniformly continuous with respect to the gradient uniformly in the $x$ variable
(in our case this amounts to assume that both controls sets are compact). Let us mention that
uniqueness and existence problems for a class of first-order Hamiltonians corresponding to
unbounded control sets and under assumptions including deterministic linear quadratic
problems have been addressed by several authors, see, e.g. the book of Bensoussan \cite{bensoussan82},
the papers of Alvarez \cite{alvarez97}, Bardi and Da Lio \cite{bdl97}, Cannarsa and Da Prato \cite{cdp89},
Rampazzo and Sartori \cite{rs00} in the case of convex operators, and the papers of
Da Lio and McEneaney \cite{dlme02} and Ishii \cite{ishii97} for more general operators.
As for second-order Hamiltonians under quadratic growth assumptions, Ito \cite{ito01} obtained
the existence of locally Lipschitz solutions to particular equations of the form \rec{hjbintr}
under more regularity conditions on the data, by establishing a priori estimates on the solutions.
Whereas Crandall and Lions in \cite{cl90} proved a uniqueness result for very particular operators
depending only on the Hessian matrix of the solution. Kobylanski \cite{kobylanski00} studied
equations with the same kind of quadratic nonlinearity in the gradient than ours, but her existence
and uniqueness results hold in the class of bounded viscosity solutions.
Finally, one can find existence and uniqueness results for viscosity solutions which may
have a quadratic growth in \cite{bbbl03} for quasilinear degenerate parabolic equations.

The second question we deal with in this paper concerns the link between the
control problem (\ref{fctval-LQintr}) and Equation (\ref{LQ-eq-intro}).
The rigorous connection between the Hamilton-Jacobi-Bellman and optimal control is usually performed
by means of a principle of optimality.
For deterministic control problems which lead to a first-order Hamilton-Jacobi
equation, see Bardi and Capuzzo Dolcetta \cite{bcd97} and Barles \cite{barles94};
for the connections between stochastic control problems and
second-order Hamilton-Jacobi-Bellmann equations we refer to
Fleming and Rishel \cite{fr75}, Krylov \cite{krylov80}, Lions \cite{lions83a, lions83b, Lions83c},
Fleming and Soner \cite{fs93}, Yong and Zhou \cite[Theorem 3.3]{yz99}
and the references therein.

However, for stochastic differential equations
with unbounded controls as in stochastic linear quadratic problems, additional difficulties arise.
Some results in this direction were obtained for infinite horizon problems, in the
deterministic case by Barles \cite{barles90b} and in the stochastic case by Alvarez \cite{alvarez96}.
In this paper, we characterize the value function (\ref{fctval-LQintr}) as the unique solution
of (\ref{LQ-eq-intro}).
Actually, our results apply for a larger class of unbounded stochastic control problems
\begin{equation}
 V(x,t)=\inf_{\alpha_s \in {\cal A}_t} \mbox{E}\{\int_t^T
\ell(X_s,s,\alpha_s)\,ds+\psi(X_T)\}
\end{equation}
where the process $X_s$ is governed by
\begin{equation}
\left\{\begin{array}{l}
dX_s=b(X_s,s,\alpha_s)ds+\sigma(X_s,s,\alpha_s)dW_s \\
X_t=x,
\end{array}\right.
\end{equation}
where $A$ is a possibly unbounded subset of a normed linear space and
all the datas are continuous with the following restricted growths:
$b$ grows at most linearly with respect to both the
control and the state, $\sigma$ grows at most linearly with respect to the state and is
bounded in the control variable, $\psi$ can have a quadratic growth and $\ell$ grows
at most quadratically with respect to both the
control and the state with a coercitivity assumption
\begin{equation}\label{coercitive-ell}
\ell (x,\alpha,t) \geq \frac{\nu}{2} |{\alpha}|^2 -C(1+|x|^2).
\end{equation}
In this case, the Hamilton-Jacobi equation looks like (\ref{hjbintr}) with $G\equiv 0$
(see Section \ref{sec-applications} for details).
Because of the unbounded framework, the use of an optimality principle to
establish the connection between the control problem and the equation is more
delicate than usual. Thus we follow another strategy which consists in comparing
directly the value function with the unique solution of the Hamilton-Jacobi equation
as long as this latter exists (see Theorem \ref{valeur-unique-sol}).

It is worth noticing that, surprisingly,
for the general stochastic linear problem,
it is not even clear how to give sense to the partial
differential equation associated to the stochastic control problem! For instance, consider
again the above linear quadratic problem where $\sigma (x,t,\alpha)= C(t)x+D(t)$ is replaced by
$C(t)x+D(t)\alpha$ which depends now both on the state and the control.
Taking $A,C,Q\equiv 0$ and $B, R,D\equiv 1$ to simplify, the Hamiltonian in (\ref{LQ-eq-intro})
becomes
\begin{equation}
\mathop{\rm sup}_{\alpha\in \R^k}
\Big\{ -\langle \alpha ,  Dw \rangle
- |\alpha|^2 -\frac{|\alpha|^2}{2}\Delta w  \Big\}
\end{equation}
which is $+\infty$ as soon as $\Delta w \leq -2.$
The connection between the control problem and the equation in this case was already
investigated (see Yong and Zhou \cite{yz99} and the references therein).
The results need {\it a priori}  knowledges about the value function (the
value function and its derivatives are supposed to remain in the domain of the Hamiltonian).
We do not consider the case when $\sigma$ is unbounded in the control
variable in this paper. It is the aim of a future work.
When finishing this paper, we learnt that Krylov \cite{krylov01} succeeded in treating the
general stochastic linear regulator. But his assumptions are designed to solve this latter
problem (the datas are supposed to be polynomials of degree 1 or degree 2 in $(x,\alpha)$)
and the proofs rely heavily on the particular form of the datas.

Another important example of equations of type (\ref{hjbintr}) where concave and convex
Hamiltonians appear is the first-order equation
\begin{equation}
\displaystyle\frac{\partial  w}{\partial t}
+\min_{\alpha\in\R^k}\left\{\frac{\gamma^2}{2}|\alpha|^2-\langle \sigma(x)\alpha ,
D w\rangle \right\}
+\max_{\beta\in B}\{-\langle g(x,\beta), D w\rangle -f(x,\beta)\}=0
\end{equation}
in $\R^N\times (0,T).$ Such kind of equation is related to the so-called $H_\infty$-Robust control problem.
This problem can be seen as a deterministic differential game.
We refer the reader to Example \ref{example-risk-sensitive} and 
McEneaney \cite{mceneaney95, mceneaney95b, mceneaney98}, Nagai \cite{nagai96}
and the references therein for details.

Finally, we point out that one of the main fields of application of these types of equations and problems
is mathematical finance, see e.g. Lamberton and Lapeyre \cite{ll97}, Fleming and Soner \cite{fs93},
{\O}ksendal \cite{oksendal98} and the references therein for an introduction.
For recent papers which deal with equations we are interested in, we refer to Pham \cite{pham02} and
Benth and Karlsen \cite{bk03} (see Example \ref{exple-math-finance}).

Let us now describe how the paper is organized.

Section \ref{sec-1-comp} is devoted to the study of (\ref{hjbintr}). More precisely,
we prove a uniqueness result for (\ref{hjbintr}) in the set of continuous
functions growing at most quadratically in the state variable under the assumption that either $A$ or
$B$ is an unbounded control set, the functions $b, g$ and $\ell, f$ grow respectively at most linearly and
quadratically with respect to both the control and the state. Instead the functions $\sigma , c$  are
assumed to grow at most linearly with respect to the state and  bounded in the control variable.

One of the main tools within the theory of viscosity solutions to obtain  a uniqueness result is to
show a comparison result between viscosity upper semicontinuous subsolutions and lower semicontinuous
supersolutions to \rec{hjbintr}, see Theorem \ref{thm-unicite}. Indeed,  the existence and the
uniqueness (Corollary \ref{thm-existence}) follow as a by-product of the comparison result and
Perron's method of Ishii \cite{ishii87}. However, under our general assumptions one cannot expect
the existence of a solution for all times as Example \ref{LQ-deterministe} shows.

The method we use in proving the comparison Theorem \ref{thm-unicite} is similar in the spirit to
the one applied  by Ishii in \cite{ishii97} in the case of first order Hamilton-Jacobi equations and it is based
on a kind of linearization
procedure of the equation. Roughly speaking, it consists in three main steps: 1) one computes the
equations satisfied by $\Psi_\mu=U-\mu V$ ($U,V$ being respectively the sub
and supersolution of the original pde and $0< \mu < 1$ a parameter); 2) for all $R>0$ one constructs
a strict supersolution
$\chi^\mu_R$ of the linearized equation such that $\chi^\mu_R(x,t)\to 0$ as $R\to\infty$; 3) one shows
that $\Psi_\mu(x,t)\le  \chi^\mu_R(x,t)$ and then one concludes by letting first $R\to\infty$ and then
$\mu\to 1.$ 

In Section \ref{sec-applications}, we give applications to finite horizon stochastic
control problems previously mentioned and we provide some examples.

In Section \ref{sec-related}, we deal with
particular cases where both controls are unbounded but $H$ (or $G$) is ``predominant" in $H+G$
(see Remark \ref{extension} and Theorem \ref{thm-related}).
For instance, we are able to deal with equations of the form
$$
\displaystyle\frac{\partial  w}{\partial t}
-\frac{|\Sigma_1(x) D w|^2}{2}+\frac{|\Sigma_2(x)
D w|^2}{2}=0~~\mbox{in $\R^N\times(0,T)$}
$$
where $\Sigma_1,~\Sigma_2$ are $N\times k$ matrices,
which corresponds to the case  $\alpha, \beta\in\R^k,$
$\sigma \equiv c\equiv 0$, $b(x,t,\alpha)= \Sigma_1(x)\alpha,$ $g(x,t,\beta)= \Sigma_2(x)\beta$ and
$\ell(x,t,\alpha)={|\alpha|^2}/{2},$ $f(x,t,\beta)={|\beta|^2}/{2}$ in (\ref{expression-H-intr})
and (\ref{expression-G-intr}). 
The comparison result applies
if either $(\Sigma_1 \Sigma_1^T)(x) > (\Sigma_2 \Sigma_2^T)(x)$ or 
$(\Sigma_1 \Sigma_1^T)(x) < (\Sigma_2 \Sigma_2^T) (x).$

When neither $H$ nor $G$) is ``predominant", the problem seems to be very difficult
and our only result takes place in dimension $N=1$: 
we have comparison for
\begin{eqnarray} \label{eikonal-quad}
\displaystyle\frac{\partial  w}{\partial t}+h(x,t)|D w|^2=0 \ \ \ {\rm in} \ \R^N\times(0,T),
\end{eqnarray}
where the function $h$ may change sign (see \textbf{(A5)} for details).
Finally, we point out that assumptions and proofs in Section \ref{sec-related}
essentially differ from those of Theorem \ref{thm-unicite} (see Remark \ref{rmq-eq-rel}).

{\bf Acknowledgments.} We are grateful to Martino Bardi and Guy Barles for valuable suggestions
during the preparation of this work. The first author is partially supported by M.I.U.R. project
``Viscosity, metric, and control theoretic methods for nonlinear partial differential equations.''
The second author is partially supported by the ACI project ``Mouvement d'interfaces avec termes
non locaux''.

\section{Comparison result for the Hamilton-Jacobi equation (\ref{hjbintr})}
\label{sec-1-comp}

In order to give precise assumptions upon Equation (\ref{hjbintr}) and (\ref{expression-H-intr}), 
(\ref{expression-G-intr}), we need to introduce some notations.
For all integers $N,M\geq 1$ we denote by ${\cal{M}}_{N, M} (\R)$
(respectively ${\cal{S}}_{N} (\R),$ ${\cal{S}}_{N}^+ (\R)$)  the set of  real  $N\times M$ matrices
(respectively real symmetric matrices, real symmetric nonnegative  $N\times N$ matrices).
All the norms which appear in the sequel are denoting by $|\cdot|.$
The standard Euclidean inner product in $\R^N$ is written $\langle\cdot,\cdot\rangle.$
We recall that a modulus of continuity $m: \R\to \R^+$ is a nondecreasing continuous function
such that $m(0)=0.$ Finally $B(x,r)= \{ y\in \R^N : |x-y|<r\}$ is the open ball of center $x$ and radius $r>0.$ 

We list the basic assumptions on $H,$ $G$ and $\psi.$
We assume that there exist positive constants $\bar{C}$ and
$\nu$ such that:

{\noindent \textbf{(A1)} (Assumptions on $H$ given by (\ref{expression-H-intr}))~:}\par
\begin{enumerate}
\item[(i)]~$A$ is a subset of a separable complete normed space.
The main point here is the possible unboundedness of $A.$
Therefore, to enlight this property in the sequel, we take
$A=\R^k$ for some $k\ge 1$ (see Remark \ref{controles-non-bornes} above);
  \item[(ii)]~${b}\in C(\R^N \times [0,T]\times \R^k ; \R^N)$ satisfying for
$x,y\in\R^N, t\in [0,T], {\alpha}\in\R^k,$
\begin{eqnarray*}
|{b}(x,t,{\alpha})-{b}(y,t,{\alpha})|& \leq &\bar{C} (1+|{\alpha}|)|x-y| \\
|{b}(x,t,{\alpha})| &\leq & \bar{C} (1+|x|+|{\alpha}|)\ ;
\end{eqnarray*}
\item[(iii)]~
${\ell}\in C(\R^N \times [0,T]\times \R^k ; \R)$ satisfying for $x\in\R^N, t\in [0,T], {\alpha}\in\R^k,$
$$
\bar{C}(1+|x|^2+|{\alpha}|^2)
\geq
{\ell}(x,t,{\alpha})\geq \frac{\nu}{2} |{\alpha}|^2 +\ell_0 (x,t,{\alpha})
\ \ {\rm with} \ \ell_0 (x,t,{\alpha} ) \geq -\bar{C} (1+|x|^2),
$$
and for every $R>0,$ there exists a modulus of continuity $m_R$
such that  for all $  x,y\in B(0,R), t\in [0,T], {\alpha}\in\R^k,$
\begin{equation}
 |{\ell}(x,t,{\alpha})-{\ell}(y,t,{\alpha})| \leq (1+|{\alpha}|^2)\, m_R(|x-y|)\ ;
\end{equation}
\item[(iv)]~
${\sigma}\in C(\R^N\times [0,T]\times \R^k ; {\cal{M}}_{N,M} (\R))$ is locally
Lipschitz continuous with respect to $x$ uniformly for
$(t, {\alpha}) \in [0,T]\times \R^k$
and satisfies { for \ every} $  x\in \R^N, \ t\in [0,T], \ {\alpha}\in \R^k,$
$$
|{\sigma}(x,t, {\alpha})| \leq  \bar{C}(1+|x|).
 $$
\end{enumerate}

{\noindent \textbf{(A2)} (Assumptions on $G$ given by(\ref{expression-G-intr}))~:} \par
\begin{enumerate}
\item[(i)]~ ${{B}}$ is a bounded subset of a normed space;
\item[(ii)]~ ${g}\in
C(\R^N\times [0,T]\times {{B}}; \R^N)$ is locally Lipschitz continuous
with respect to $x$ uniformly for $(t, \beta)\in [0,T]\times {{B}}$
and satisfies  for  every $x\in \R^N,  \ t\in [0,T],\ \beta\in
{{B}},$
$$
|{g}(x,t, \beta)|\le \bar{C}(1+|x|)\ ;
 $$
\item[(iii)]~
${f}\in C(\R^N\times [0,T]\times {{B}}; \R)$ is locally uniformly continuous with respect to $x$ uniformly in
 $(t,\beta)\in [0,T]\times {{B}}$ and satisfies   for every $ x\in \R^N, \ t\in [0,T],\ \beta\in {{B}},$
$$
|{f} (x,t, \beta)| \le \bar{C}(1+|x|^2);
 $$
\item[(iv)]~ ${c}\in C(\R^N\times [0,T]\times {{B}};
{\cal{M}}_{N,M} (\R))$ is locally Lipschitz continuous with
respect to $x$ uniformly for $(t, \beta) \in [0,T]\times {{B}}$ and
satisfies   for every $ x\in \R^N, \ t\in [0,T],\ \beta\in {{B}},$
$$
|{c}(x,t, \beta)| \leq \bar{C}(1+|x|).
$$
\end{enumerate}

{\noindent \textbf{(A3)} (Assumptions on the initial condition $\psi$)~:} \\
$\psi \in C(\R^N;\R)$ and
$$
|\psi(x)| \leq \bar{C}(1+|x|^2) \ \ \ \
{\rm for \ every} \ x\in \R^N.
$$

\begin{Remark} \rm \label{controles-non-bornes}
(i) About \textbf{(A1)}(i): we choose to take $A=\R^k$ in this section to enlight
the possible unboundedness of $A$ in the notation. Indeed, the calculations
when $A$ is any subset of a complete separable normed space are the same
and are based on the following inequality:
for every $\rho>0, \gamma\in\R,$
\begin{eqnarray} \label{ineg-fondamentale}
\mathop{\rm inf}_{{\alpha}\in A}\left\{
\rho|{\alpha}|^2 + \gamma|{\alpha}|
\right\}
= \mathop{\rm inf}_{{\alpha}\in A}\left\{
\left( \sqrt{\rho} |\alpha|^2 + \frac{\gamma}{2\sqrt{\rho}}\right)^2 -\frac{\gamma^2}{4\rho}
\right\}
\geq -\frac{\gamma^2}{4\rho}.
\end{eqnarray}

(ii) Note that, with respect to the gradient variable,  $H$ is a concave function and $G$
is a convex function. Under Assumptions \textbf{(A1)} and \textbf{(A2)},
classical computations show that $H$ and $G$ are continuous in all their variables.
\end{Remark}

\begin{Example} \rm
The typical case we have in mind is
when $H$ is quadratic in the gradient variable, for instance
$A=\R^k,$ ${\ell}(x,t,{\alpha})= |{\alpha}|^2/2,$ ${\sigma}\equiv 0$ and
${b}(x,t,{\alpha}) = a(x){\alpha}$ where
$a \in C(\R^N ; {\cal{M}}_{N,k} (\R))$
is locally Lipschitz
continuous and bounded for all $x\in \R^N.$
It leads to
\begin{equation}\label{examplers}
H(x,p)= \mathop{\rm inf}_{{\alpha}\in\R^k}
\{ \langle a(x){\alpha}, p \rangle + \frac{|{\alpha} |^2}{2} \}
= -\frac{|a(x)^T p|^2}{2}.
\end{equation}
This particular example is treated both  in
 Ishii \cite{ishii97} and in Da Lio and McEneaney \cite{dlme02} in the case of  first order
Hamilton-Jacobi-Bellman equations under more restrictive assumptions than ours. In particular,
$a$ has to be a nonsingular matrix in \cite{ishii97}. See also
Section \ref{sec-related} for some further comments.
\end{Example}

For any $O\subseteq\R^K$, we denote by $USC(O)$ the set of
upper semicontinuous functions in $O$ and by $LSC(O)$ the set of
lower semicontinuous functions in $O.$
 \par The main result of this section is
the

\begin{Theorem} \label{thm-unicite}
Assume \textbf{(A1)}--\textbf{(A3)}.
Let $U \in USC(\R^N\times [0,T])$ be a viscosity subsolution of (\ref{hjbintr})
and $V \in LSC(\R^N\times [0,T])$ be a viscosity supersolution of (\ref{hjbintr}).
Suppose that $U$ and $V$ have quadratic growth, i.e.
there exists $\hat{C}>0$ such that,
for all $x\in \R^N,$ $t\in [0,T],$
\begin{eqnarray}\label{croiss-quad}
|U(x,t)|, \, |V(x,t)| \, \leq \hat{C}(1+|x|^2).
\end{eqnarray}
Then $U\leq V$ in $\R^N\times [0,T].$
\end{Theorem}

The question of the existence of a continuous solution to  (\ref{hjbintr}) is
not completely obvious. In the framework of viscosity solutions, existence
is usually obtained as a consequence of the comparison principle by
means of Perron's method as soon as we can build a sub- and a super-solution
to the problem.
Here, the comparison principle is proved in the class of functions satisfying
the quadratic growth condition (\ref{croiss-quad}).
Therefore, to perform the above program of existence, we need to
be able to build quadratic sub- and super-solutions to (\ref{hjbintr}).
In general one can expect to build such sub- and super-solutions only
for short time (see the following lemma and Corollary \ref{thm-existence}).
In Example \ref{LQ-deterministe}, we see that solutions may not exist
for all time.
\begin{Lemma} \label{lemme-existence}
Assume \textbf{(A1)}--\textbf{(A3)}. If $K\geq \bar{C}+1$ and $\rho$ are
large enough, then $\underline u(x,t)=-K{\rm e}^{\rho t}(1+|x|^2)$ is a viscosity
subsolution of (\ref{hjbintr}) in $\R^N\times [0,T]$  and
there exists $0< \tau \leq T$ such that
$\overline u(x,t)=K{\rm e}^{\rho t}(1+|x|^2)$ is a viscosity
supersolution of (\ref{hjbintr}) in $\R^N\times [0,\tau].$
\end{Lemma}

{\it Proof of Lemma \ref{lemme-existence}.}
We  only verify that $\overline u$ is a supersolution (the proof that
$\underline u$ is a subsolution being similar and simpler).
Since $K\geq \bar{C}+1$, we have
$\overline u(x,0)=K(1+|x|^2)\ge \psi(x).$
Moreover, since $\overline{u}$ is smooth and using \textbf{(A1)},
\textbf{(A2)} and (\ref{ineg-fondamentale}), we have
\begin{eqnarray*}
& & \displaystyle\frac{\partial \overline u}{\partial t}
+ H(x,t,D\overline u , D^2 \overline u)
+ G(x,t,D\overline u, D^2 \overline u) \\
& = &
K\rho {\rm e}^{\rho t}(1+|x|^2)
+ H(x,t,2K {\rm e}^{\rho t} x, 2K {\rm e}^{\rho t} Id)
+ G(x,t,2K {\rm e}^{\rho t} x, 2K {\rm e}^{\rho t} Id) \\
& \geq &\displaystyle
K\rho {\rm e}^{\rho t}(1+|x|^2)
+\mathop{\rm inf}_{{\alpha}\in \R^k}
\{ -\bar{C}(1+|x|+|{\alpha}|)\}2K {\rm e}^{\rho t} |x|
+\frac{\nu}{2}|{\alpha}|^2- \bar{C}(1+|x|^2)
- \bar{C}^2 (1+|x|)^2 2K {\rm e}^{\rho t} \} \\
& & \hspace*{3.5cm}
+ \mathop{\rm sup}_{\beta\in {{B}}}
\{ -\bar{C}(1+|x|)2K {\rm e}^{\rho t} |x|- \bar{C}(1+|x|^2)
- \bar{C}^2 (1+|x|)^2 2K {\rm e}^{\rho t} \} \\
& \geq &\displaystyle
K\rho {\rm e}^{\rho t}(1+|x|^2) -K(10\bar{C}+12\bar{C}^2){\rm e}^{\rho t}(1+|x|^2)
+ \mathop{\rm inf}_{{\alpha}\in \R^k}
\{ -2K \bar{C}{\rm e}^{\rho t} |x||{\alpha}|+\frac{\nu}{2}|{\alpha}|^2 \} \\
& \geq &\displaystyle
\left[\rho - 10\bar{C}-12\bar{C}^2-\frac{2\bar{C}^2 K {\rm e}^{\rho t}}{\nu}\right]
K {\rm e}^{\rho t} (1+|x|^2).
\end{eqnarray*}
We notice that if $t\rho\leq 1$ and $\rho>0$ is large enough, then
the quantity between the brackets is nonnegative. Hence the result
follows with $0< \tau=1/ \rho.$~~$\Box$

As explained above, Theorem \ref{thm-unicite} together with
Perron's method implies the following result.
We omit its proof since it is standard.
\begin{Corollary} \label{thm-existence}
Assume \textbf{(A1)}--\textbf{(A3)}. Then there is $\tau >0$ such that
there exists a unique continuous viscosity solution
of (\ref{hjbintr}) in $\R^N\times [0,\tau ]$ satisfying the growth
condition  (\ref{croiss-quad}).
\end{Corollary}

\begin{Remark} \rm \label{extension} 
(i) For global existence results under further
regularity assumptions on the data
we refer the reader to Ito \cite{ito01}.
For a case where blowup in finite time occurs, see
Example \ref{LQ-deterministe} which follows. 

(ii) Theorem \ref{thm-unicite} and Corollary \ref {thm-existence}
hold when replacing ``inf'' by ``sup'' and ``${\ell}$'' by ``$-{\ell}$''
in $H$ or/and ``sup'' by ``inf'' and ``${f}$'' by ``$- {f}$''
in $G.$ To adapt the proofs, one can use the
change of function $ w^\prime :=- w$ or to consider
$\mu U -V$ instead of $U-\mu V$ in the proof of Theorem  \ref{thm-unicite}.
Therefore it is possible to deal either with unbounded controls
in the ``sup'' in order to have a convex quadratic Hamiltonian
or to deal with unbounded control in the ``inf'' in order to have
a concave quadratic Hamiltonian. 

(iii) Up to replace $t$ by $T-t,$ all our results hold for
\begin{equation}
\left\{\begin{array}{ll}
\displaystyle -\frac{\partial  w}{\partial t}+
H (x,t,D w, D^2  w)
+ G (x,t,D w, D^2  w) = 0 & \mbox{in  $\R^N\times (0,T),$}\\[3mm]
 w (x,T)= \psi (x)& \mbox{in  $\R^N$}\end{array}\right.\
\end{equation}
This latter equation with terminal condition is the one which
arises usually in control theory, see Example \ref{LQ-deterministe}
and Section \ref{sec-applications}. 

(iv) In this section, we are not able to consider the case
when both ${\alpha}$ in $H$ and $\beta$ in $G$ are unbounded controls.
Roughly speaking, one of the reason is that unbounded controls lead to 
quadratic Hamiltonians. When both controls are unbounded, we then obtain
two quadratic-type Hamiltonians, a concave and a convex one.
Let us explain the difficulty on a model case where
\begin{eqnarray*}
H(x,p)= \mathop{\rm inf}_{{\alpha}\in\R^k}
\{ \langle a_1 (x){\alpha}, p \rangle + \frac{|{\alpha} |^2}{2} \}
= -\frac{|a_1 (x)^T p|^2}{2}, \
G(x,p)= \mathop{\rm sup}_{\beta\in\R^k}
\{ \langle a_2 (x)\beta, p \rangle - \frac{|\beta |^2}{2} \}
= \frac{|a_2 (x)^T p|^2}{2}, 
\end{eqnarray*}
with $a_1, a_2\in C(\R^N ; {\cal{M}}_{N,k}(\R))$ are locally
Lipschitz continuous and bounded.
The difficulty to treat such a case is related to our strategy of proof
which relies on a kind of ``linearization procedure''  
(see Lemma \ref{linearisation} and its proof). In this simple case, this ``linearization''
uses in a crucial way the convex inequality
$$
\frac{|p|^2}{\mu}-|q|^2\geq -\frac{|p-q|^2}{1-\mu} 
\ \ \ {\rm for \ all} \ p,q\in \R^N \ {\rm and} \ 0<\mu <1,
$$
which does not work at the same time for a concave and a convex Hamiltonian.
Of course, in this simple case, there are alternative ways to
solve the problem: we have
\begin{eqnarray} \label{hamiltonien-related-1}
H(x,p)+G(x,p) &=& \frac{1}{2}
\langle (a_1 +a_2 )(a_2-a_1)^T p,p\rangle,
\end{eqnarray}
which allows to apply Theorem \ref{thm-unicite} up to add
some assumptions on $a_1$ or $a_2$ (for instance $(a_1 +a_2 )(a_2-a_1)^T$
is a nonnegative symmetric matrix with a locally Lipschitz 
squareroot).
In Section \ref{sec-related}, we provide another approach to solve such
equations (see Theorem \ref{thm-related} and Remark \ref{rmq-eq-rel}).
\end{Remark}

\begin{Example} \rm \label{LQ-deterministe}
A deterministic linear quadratic control problem. \\
Linear quadratic control problems (see also Section
\ref{sec-applications} for stochastic linear quadratic control
problems) are the typical examples we have in mind since they lead
to Hamilton-Jacobi equations with quadratic terms. On the other hand,
the value function can blow up in finite time.
Consider the control problem (in dimension 1 for sake of simplicity)
\begin{eqnarray*}
\left\{
\begin{array}{cc}
dX_s = {\alpha}_s\, ds, & s\in [t,T], \ 0\leq t\leq T, \\
X_t = x\in \R, &
\end{array}
\right.
\end{eqnarray*}
where the control ${\alpha} \in {\cal{A}}_t := L^2 ([t,T] ; \R)$
and the value function is given by
\begin{eqnarray*}
V(x,t)= \mathop{\rm inf}_{{\alpha} \in {\cal{A}}_t}
\{ \rho \int_t^T (|{\alpha}_s|^2 + |X_s|^2)\, ds - |X_T|^2\} \ \
\ {\rm for \ some} \ \rho >0.
\end{eqnarray*}
Then the value function, when it is finite, is the unique viscosity solution
of the Hamilton-Jacobi equation of type (\ref{hjbintr}) which reads
\begin{eqnarray} \label{hjb1-lq}
\left\{
\begin{array}{cc}
-  w_t +\frac{1}{4\rho} | w_x|^2 = \rho x^2 \ \
& {\rm in} \ \R\times (0,T), \\
 w (x,T)= -x^2, &
\end{array}
\right.
\end{eqnarray}
(see Theorem \ref{valeur-unique-sol} for a proof of this result).
Looking for a solution
$ w$ under the form $ w (x,t)= \varphi (t)x^2,$ we obtain that
$\varphi$ is a solution of the differential equation
\begin{eqnarray}\label{edo-control}
-\varphi ' + \frac{\varphi^2}{\rho}=\rho \ \  {\rm in} \ (0,T), \ \ \ \ \
\varphi (T)= -1.
\end{eqnarray}
We distinguish two cases: \\
Case 1.--- If $\rho \geq 1,$ then the solution of (\ref{edo-control}) is defined
in the whole interval $[0,T]$ for all $T>0$ and is given by
\begin{eqnarray}\label{sol-edo}
\varphi (t)= \rho \,
\frac{(\rho -1){\rm e}^{2(T-t)} - (\rho +1)}
{(\rho -1){\rm e}^{2(T-t)} + \rho +1}.
\end{eqnarray}
which is a function decreasing from $\varphi (0)$ to $-1.$
Therefore (\ref{hjb1-lq}) admits a unique viscosity solution
in $\R\times [0,T]$ which is the value function of the control
problem
\begin{eqnarray*}
V(x,t)= \varphi(t)x^2.
\end{eqnarray*}
Note that, if $\rho >1$ and $T> {\rm ln}((\rho +1)/(\rho -1))/2,$ then $\varphi (0)>0.$ It follows
that the value function satisfies (\ref{croiss-quad}) but
is neither bounded from above neither bounded from below. \\
Case 2.--- If $0< \rho < 1$ and $T> {\rm ln}((1+\rho )/(1-\rho ))/2,$
then the solution of (\ref{edo-control}) is given by (\ref{sol-edo}) in
$( \bar\tau ,T],$ where
$$
\bar\tau := T- \frac{1}{2}\, {\rm ln} \left( \frac{1+\rho}{1-\rho}\right),
$$
and blows up at $t=\bar\tau.$
Therefore we have existence for (\ref{hjb1-lq}) only in $\R\times (\bar\tau, T].$
\end{Example}


{\it Proof of Theorem \ref{thm-unicite}.}
We divide the proof of the theorem in two steps.

{\it Step 1.}  We first assume that $\ell_0 \geq 0,$
${f} \leq 0$ and $\psi \leq 0$ in \textbf{(A1)}, \textbf{(A2)}
and \textbf{(A3)}. \\
The proof is based on the two following lemmas whose proofs
are postponed.
\begin{Lemma} \label{linearisation}
{(``Linearization'' of Equation (\ref{hjbintr})).}
Let $0<\mu <1$ and set $\Psi=U-\mu V.$ Suppose
$\ell_0 \geq 0,$ ${f} \leq 0$ and $\psi \leq 0.$ Then $\Psi$
is an $USC$ viscosity subsolution of
\begin{eqnarray}\label{hjb2}
{\cal{L}}[ w ]  & := &  \
\frac{\partial  w}{\partial t}-
\frac{\bar{C}^2}{2\nu (1-\mu)}
| D w |^2
- 2 \bar{C}(1+|x|) |D w|
\\[2mm]
& & - \mathop{\rm sup}_{{\alpha}\in \R^k} {\rm Tr}\left[ {\sigma}(x,t,
{\alpha}){\sigma}^T(x,t,{\alpha}) D^2  w \right] - \mathop{\rm sup}_{\beta\in
{{B}}} {\rm Tr}\left[ {c}(x,t, \beta){c}^T(x,t,\beta) D^2
 w \right] \ = \ 0 \nonumber
\end{eqnarray}
in $\R^N\times (0,T),$ with the initial condition
\begin{equation} \label{hjb2-initial}
 w (\cdot,0)\leq (1-\mu)\psi \leq 0.
\end{equation}
\end{Lemma}

\begin{Lemma} \label{prob-parab}
Consider the parabolic problem
\begin{equation}\label{equa-parab}
\left\{\begin{array}{ll}
 \varphi_t -r^2 \varphi_{rr} - r
\varphi_{r} =0
 & {  \mbox{in}}\ \
[0,+\infty)\times (0,T]\\
\varphi (r,0)= \varphi_R (r) &  \mbox{in}\ \ [0,+\infty),
\end{array}\right.
\end{equation}
where $\varphi_R (r) = {\rm max} \{0,r-R\}$ for some $R>0.$
Then (\ref{equa-parab}) has a unique solution
$\varphi\in C([0,+\infty)\times [0,T]) \cap
C^\infty ([0,+\infty)\times (0,T])$ such that,
for all $t\in (0,T],$ $\varphi (\cdot,t)$ is positive, nondecreasing
and convex in $[0,+\infty).$ Moreover, for every $(r,t)\in
[0,+\infty)\times (0,T],$
\begin{eqnarray} \label{derivee-phi}
\varphi (r,t)\geq  \varphi_R (r), \ \ \
0\leq \varphi_r (r,t)\leq {\rm e}^T \ \ \ { and} \ \ \
\varphi (r,t) \mathop{\longrightarrow}_{R\to +\infty} 0.
\end{eqnarray}
\end{Lemma}

Let $\Phi (x,t)= \varphi (C(1+|x|^2) {\rm e}^{Lt}, Mt)+\eta t$
where $\varphi$ is given by Lemma \ref{prob-parab}, $L,M, \eta$ are
positive constants to be determined and
$C > {\rm max}\{ \bar{C}, \hat{C}\},$
where $\bar{C}$ and $\hat{C}$ are the constants appearing
in the assumptions.

{\it Claim~:}{\it We can choose the constants $L$ and $M$ such that
$\Phi$ is a strict supersolution
of (\ref{hjb2}) at least in $\R^N\times (0,\tau]$ for small $\tau.$}

To prove the claim, we have to show that ${\cal{L}}[\Phi]>0$
in $\R^N\times  (0,\tau]$ for some $\tau >0.$
The function $\Phi \in C(\R^N\times [0,T]) \cap
C^\infty (\R^N\times (0,T])$ and we have, for $t>0,$
\begin{eqnarray*}
\Phi_t = \eta+LC(1+|x|^2){\rm e}^{Lt}\varphi_r +M\varphi_t,
\ \ 
D\Phi = 2Cx{\rm e}^{Lt} \varphi_r \ \ 
{\rm and}  \ \ D^2 \Phi = 2C \, Id \, {\rm e}^{Lt}\varphi_r
+ 4C^2 {\rm e}^{2Lt} \varphi_{rr} \, x\otimes x .
\end{eqnarray*}
Using \textbf{(A1)} and  \textbf{(A2)},
for all $(x,t)\in \R^N\times (0,T],$ we get
\begin{eqnarray*}
{\cal{L}}[\Phi ] & \geq &
\eta+LC(1+|x|^2){\rm e}^{Lt}\varphi_r + M\varphi_t
- \frac{\bar{C}^2}{2\nu(1-\mu)}
| 2C x {\rm e}^{Lt} \varphi_r |^2
- 2 \bar{C}(1+|x|)|2Cx{\rm e}^{Lt} \varphi_r| \\
& &
-2 \bar{C}^2 (1+|x|)^2  |2C \, Id \, {\rm e}^{Lt}\varphi_r
+ 4C^2 {\rm e}^{2Lt} x\otimes x \, \varphi_{rr} |.
\end{eqnarray*}
Setting $r= C(1+|x|^2){\rm e}^{Lt}$ and since $C>\bar{C},$
we obtain
\begin{eqnarray*}
{\cal{L}}[\Phi ] & \geq & \eta+
M \varphi_t - 16 C^2 r^2 \varphi_{rr} - 8C(C+1) r \varphi_r
+ \left( L- \frac{2C^3{\rm e}^{Lt} \varphi_r }{\nu(1-\mu)}
\right)  r \varphi_r.
\end{eqnarray*}
Our aim is to fix the parameters $M$ and $L$ in order to make
${\cal{L}}[\Phi ]$ positive.

We first choose $M > 16C^2+  8C.$ Since $\varphi$ is a solution
of (\ref{equa-parab}) (Lemma \ref{prob-parab}), we obtain
\begin{equation}
{\cal{L}}[\Phi ] > \eta +
\left( L- \frac{2C^3{\rm e}^{Lt} \varphi_r }{\nu(1-\mu)}\right)
r \varphi_r.
\end{equation}
Then, taking $\displaystyle
L>  \frac{2C^3{\rm e}^{T+1}}{\nu(1-\mu)},
$
we get ${\cal{L}}[\Phi ] > \eta >0$ for all $x\in\R^N$ and $t\in (0,\tau],$
where $\tau = 1/L.$ This proves the claim.

We continue by considering
\begin{eqnarray}\label{max-sous-sur}
\mathop{\rm max}_{\R^N\times [0,\tau]} \{ \Psi-\Phi \},
\end{eqnarray}
where $\Psi$ is the function defined in Lemma \ref{linearisation}
which is a  viscosity subsolution of (\ref{hjb2})
and $\Phi$ is the strict supersolution of
(\ref{hjb2}) in $\R^N\times (0,\tau]$ we built above.

From (\ref{derivee-phi}), we have
$\Phi (x,t) \geq C(1+|x|^2) > \bar{C}(1+|x|^2) \geq \Psi(x,t)$
for $|x|\geq R.$ It follows that the maximum (\ref{max-sous-sur})
is achieved at a point $(\bar{x},\bar{t})\in \R^N\times [0,\tau].$
We claim that  $\bar{t}=0$. Indeed suppose by contradiction that $\bar{t}>0$. Then  since $\Psi$
is a viscosity subsolution of  (\ref{hjb2}), by taking $\Phi$ as a
test-function, we would have
${\cal{L}}[\Phi ](\bar{x},\bar{t}) \leq 0$ which contradicts
the fact that $\Phi$ is a strict supersolution.

Thus, for all $(x,t)\in \R^N\times [0,\tau],$
$$
\Psi(x,t)-\Phi (x,t)\leq \Psi(\bar{x},0)- \Phi (\bar{x},0)
\leq (1-\mu )\psi (\bar{x}),
$$
where the last inequality follows from (\ref{hjb2-initial}) and the fact that $\Phi \geq 0.$
Since we assumed that $\psi$ is nonpositive,  for every $(x,t)\in \R^N\times [0,\tau],$ we have
$
\Psi(x,t)\leq \Phi (x,t).
$
Letting $\eta$ go to $0$ and $R$ to $+\infty,$ we get by Lemma \ref{prob-parab},
$\Psi\leq 0$ in $\R^N\times [0,\tau].$

By a step-by-step argument, we prove that
$\Psi\leq 0$ in $\R^N\times [0,T].$
Therefore $\Psi= U-\mu V \leq 0$ in $\R^N\times [0,T].$
Letting $\mu$ go to 1, we obtain $U\leq V$ as well which concludes the
Step 1. 

{\it Step 2.} The general case. \\
The idea is  to reduce to  the first case by a suitable change of
function (see Ishii \cite{ishii97}).
Suppose that $ w$ is a solution of (\ref{hjbintr}). Then,
a straightforward computation shows that
$\bar{ w} (x,t)=  w (x,t) - C(1+|x|^2){\rm e}^{\rho t},$
for $C> \bar{C}, \hat{C}$ and $\rho >0,$
is a solution of
\begin{equation}
\begin{cases}
\displaystyle{
\bar{ w}_t
+ \mathop{\rm inf}_{{\alpha}\in \R^k}
\left\{
\langle {b}(x,t, {\alpha} ), D\bar{ w} \rangle
+ \bar{{\ell}}(x,t,{\alpha})
- {\rm Tr}\left[ {\sigma} (x,t,{\alpha}){\sigma}^T(x,t,{\alpha})
D^2 \bar{ w} \right]
\right\} } \\
\displaystyle{ \hspace*{2cm} + \mathop{\rm sup}_{\beta\in  {{B}}}
\left\{ -\langle {g}(x,t,\beta ), D\bar{ w} \rangle -
\bar{{f}} (x,t,\beta) -{\rm Tr}\left[ {c} (x,t,\beta)
{c}^T(x,t,\beta) D^2 \bar{ w} \right] \right\} =0,
} \\
\bar{ w} (x,0)= \psi (x)- C(1+|x|^2),
\end{cases}
\end{equation}
where
\begin{eqnarray*}
\bar{{\ell}}(x,t,{\alpha} ) = 
{\ell}(x,t,{\alpha} ) + 2C{\rm e}^{\rho t} \langle {b}(x,t,{\alpha} ), x\rangle
- 2C{\rm e}^{\rho t}\,
{\rm Tr}\left[{\sigma} (x,t,{\alpha}) {\sigma}^T(x,t,{\alpha})\right] 
+\displaystyle \frac{1}{2} C\rho
{\rm e}^{\rho t} (1+|x|^2), \\
\bar{{f}}(x,t,\beta ) = 
{{f}} (x,t,\beta) + 2C {\rm e}^{\rho t}
\langle {g}(x,t,\beta ), x\rangle + 2C{\rm e}^{\rho t}\,
{\rm Tr}\left[{c} (x,t,\beta){c}^T(x,t,\beta)\right] 
-\displaystyle  \frac{1}{2} C \rho {\rm e}^{\rho t}(1+|x|^2).
\end{eqnarray*}
We observe that   $\bar{{\ell}}$ and $\bar{{f}}$ still satisfy
respectively the assumptions \textbf{(A1)}(iii) and \textbf{(A2)}(iii).
Moreover   from \textbf{(A3)}, we can
choose
 $C>\bar{C} $
in order that $\bar{\psi}\leq 0.$

Next we show that if  $\rho>0$ is chosen in a suitable way then
$\bar{{\ell}}(x,t,{\alpha})\geq \bar{\nu} |{\alpha}|^2 /2$  for all $(x,t,{\alpha})
\in \R^N\times [0,1/\rho] \times \R^k$   and ${f}(x,t,\beta)\le
0$  for all $(x,t,\beta) \in \R^N\times [0,1/\rho] \times {{B}}.$
Indeed
  for all $(x,t,{\alpha}) \in \R^N\times [0,1/\rho]\times \R^k$ by \textbf{(A1)}
 we have
\begin{eqnarray*}
\bar{{\ell}}(x,t,{\alpha}) &\geq & \frac{\nu}{2}|{\alpha}|^2 -\bar{C}(1+|x|^2)
- 2C\bar{C} {\rm e}^{\rho t} (1+|x|+|{\alpha}|) |x|
- 2C\bar{C}^2 {\rm e}^{\rho t} (1+|x|)^2  +
\frac{1}{2} C \rho {\rm e}^{\rho t}(1+|x|^2) \\
& \geq &
\frac{\nu}{2}|{\alpha}|^2 -(\bar{C}+4C\bar{C}+4C\bar{C}^2){\rm e}^{\rho t}
(1+|x|^2)
- 2C\bar{C} {\rm e}^{\rho t}  |{\alpha}||x|
+ \frac{1}{2} C \rho {\rm e}^{\rho t}(1+|x|^2).
\end{eqnarray*}
But
$$
2C\bar{C} {\rm e}^{\rho t}  |{\alpha}||x| \leq \frac{\nu}{4}|{\alpha}|^2
+\frac{16C^2\bar{C}^2 {\rm e}^{2\rho t} }{\nu}|x|^2\, ;
$$
Therefore,  by choosing
\begin{eqnarray} \label{ro1}
\rho > 2\frac{\bar{C}}{ C} + 8\bar{C}+8\bar{C}^2+
\frac{32C\bar{C}^2 {\rm e} }{\nu},
\end{eqnarray}
we have
$$
\bar{{\ell}}(x,t,{\alpha}) \geq \frac{\nu}{4}|{\alpha}|^2
\ \ \ {\rm for \ all} \ (x,t,{\alpha})\in \R^N\times [0,1/\rho]
\times \R^k,
$$
which is the desired estimate with $\bar{\nu}= \nu/2 >0$ and
$\ell_0\equiv 0.$

The next step consists in choosing $\rho$ such that
$\bar{{f}}\le 0$ .
Using \textbf{(A2)}, the same kind of calculation as above
shows that, taking
\begin{eqnarray} \label{ro2}
\rho > 8(\bar{C}+ \bar{C}^2)+2
\end{eqnarray}
ensures $\bar{{f}}\leq 0.$

Finally if we choose  $C>\bar C$   and $\rho$ as the maximum of
the two quantities appearing in (\ref{ro1}) and (\ref{ro2}), we
are in the framework of Step 1 in $\R^N\times [0,\rho].$ Setting
$\bar{U}= U- C(1+|x|^2){\rm e}^{\rho t}$ and $\bar{V}= V-
C(1+|x|^2){\rm e}^{\rho t},$  from Step 1 we get $\bar{U}\leq
\bar{V}$ in $\R^N\times [0,1/\rho ]$; thus $U\leq V$ in
$\R^N\times [0,1/\rho ].$ Then by a step-by-step  argument  we
obtain the comparison in $\R^N\times [0,T].$ ~~$\Box$\par

\begin{Remark}\rm
A key fact in the proof to build a strict supersolution of
(\ref{hjb2}) is to use a function which is the solution of
the auxiliary -- and simpler -- pde (\ref{equa-parab}).
This idea come from mathematical finance to deal with
equations related to Blake and Scholes formula.
See for instance Lamberton and Lapeyre \cite{ll97}
and Barles {\it et al.} \cite{bbrs95}.
\end{Remark}

We turn to the proof of the lemmas  \ref{linearisation} and \ref{prob-parab}. 

{\it Proof of Lemma
\ref{linearisation}.} For $0< \mu <1,$ let $\tilde{V}=\mu V$ and
$\Psi=U-\tilde{V}.$
We divide the proof in different steps. 

{\it Step 1.} A new equation for $\tilde{V}.$ \\
It is not difficult to see that,
if $V$ is a supersolution of (\ref{hjbintr}), then $\tilde{V}$ is a
supersolution of
\begin{equation}
\begin{array}{ll}
\displaystyle{\tilde{V}_t
+\mu\, H\big(x,t, \frac{D\tilde{V}}{\mu}, \frac{D^2\tilde{V}}{\mu}\big)
+\mu\, G\big(x,t, \frac{D\tilde{V}}{\mu}, \frac{D^2\tilde{V}}{\mu}\big)
\geq 0}
 & {\rm in} \ \R^N\times (0,T)
 \\
\tilde{V} (x,0)\geq  \mu\, \psi (x)
& {\rm in} \ \R^N.
\end{array}
\end{equation}

{\it Step 2.} Viscosity inequalities for $U$ and $\tilde{V}.$ \\
This step is classical in viscosity theory.
Let $\varphi\in C^2 (\R^N\times (0,T])$ and
$(\bar{x},\bar{t})\in \R^N\times (0,T]$ be a local maximum
of $\Psi -\varphi.$ We can assume that this maximum is strict
in the same ball $\overline{B}(\bar{x},r)\times [\bar{t}-r,\bar{t}+r]$
(see \cite{barles94} or \cite{bcd97}). Let
$$
\Theta (x,y,t) = \varphi (x,t)
+\frac{|x-y|^2}{\varepsilon^2}
$$
and consider
$$
M_\varepsilon := \mathop{\rm max}_{x,y\in \overline{B}(\bar{x},r), \,
t\in [\bar{t}-r,\bar{t}+r]}
\{ U(x,t)-\tilde{V} (y,t)-\Theta (x,y,t) \}.
$$
This maximum is achieved at a point $(x_\varepsilon,
y_\varepsilon, t_\varepsilon)$ and, since the maximum is strict,
we know (\cite{barles94}, \cite{bcd97}) that
$$
\frac{|x_\varepsilon -y_\varepsilon|^2}{\varepsilon^2}\to 0 \ \ \
{\rm as} \ \varepsilon\to 0,
$$
and
$$
M_\varepsilon = U(x_\varepsilon,t_\varepsilon) -\tilde{V}
(y_\varepsilon,t_\varepsilon) -\Theta
(x_\varepsilon,y_\varepsilon,t_\varepsilon) \
\mathop{\longrightarrow}  \ U (\bar{x},\bar{t})-\tilde{V}
(\bar{x},\bar{t})-\varphi (\bar{x},\bar{t}) = \Psi
(\bar{x},\bar{t})-\varphi (\bar{x},\bar{t}).
$$
It means that, at the limit $\varepsilon\to 0,$ we obtain some
information on $\Psi -\varphi$ at $(\bar{x},\bar{t})$ which will
provide the new equation for $\Psi.$ Before that, we can take
$\Theta$ as a test-function to use the fact that $U$ is a subsolution
and $\tilde{V}$ a supersolution.
Indeed $(x,t)\in \overline{B}(\bar{x},r)\times [\bar{t}-r,\bar{t}+r]
\mapsto U(x,t)-\tilde{V} (y_\varepsilon,t)-\Theta (x,y_\varepsilon,t)$
achieves its maximum at $(x_\varepsilon,t_\varepsilon)$
and $(y,t)\in \overline{B}(\bar{x},r)\times  [\bar{t}-r,\bar{t}+r]
\mapsto -U (x_\varepsilon,t)+\tilde{V} (y,t)
+\Theta (x_\varepsilon,y,t)$
achieves its minimum at $(y_\varepsilon,t_\varepsilon).$
Thus, by  Theorem 8.3 in the User's guide \cite{cil92}  , for
every $\rho >0,$ there exist
$a_1, a_2 \in\R$ and $X,Y\in{\mathcal{S}}_{N}$ such that
\begin{eqnarray*}
\left(a_1,D_x\Theta(x_\varepsilon,y_\varepsilon,t_\varepsilon),
X\right)
  \in\bar{\mathcal{P}}^{2,+}(U)(x_\varepsilon, t_\varepsilon), 
\left(a_2,-D_y \Theta(x_\varepsilon,y_\varepsilon,t_\varepsilon), \ \ \ \
Y\right)
  \in\bar{\mathcal{P}}^{2,-}(\tilde{V})(y_\varepsilon,t_\varepsilon),
\end{eqnarray*}
$a_1-a_2=\Theta_t (x_\varepsilon,y_\varepsilon,t_\varepsilon)= \varphi_t (x_\varepsilon,t_\varepsilon)$ and
\begin{equation}\label{Ineq-Matrice}
-(\frac{1}{\rho}+|{{M}}|)I \leq
\left(
\begin{array}{cc}
X & 0\\
0 & -Y\\
\end{array}
\right)\leq {{M}}+ \rho {{M}}^2
\ \ \ {\rm where} \
{{M}}=D^2\Theta (x_\varepsilon,y_\varepsilon,t_\varepsilon ).
\end{equation}
Setting $\displaystyle{p_\varepsilon =
2\frac{x_\varepsilon-y_\varepsilon}{\varepsilon^2}},$  we have
$$
D_x\Theta(x_\varepsilon,y_\varepsilon,t_\varepsilon)
=  p_\varepsilon + D \varphi (x_\varepsilon,t_\varepsilon)
\ \ \ {\rm and} \ \ \
 D_y\Theta(x_\varepsilon,y_\varepsilon,t_\varepsilon)
= -p_\varepsilon,
$$
and
$$
{{M}}= \left(
\begin{array}{cc}
D^2 \varphi (x_\varepsilon,t_\varepsilon)+2I/\varepsilon^2
& -2I/\varepsilon^2\\
-2I/\varepsilon^2 & 2I/\varepsilon^2\\
\end{array}
\right).
$$
Thus, from (\ref{Ineq-Matrice}), it follows
\begin{eqnarray} \label{ineq-fq}
\langle X p,p\rangle - \langle Y q,q\rangle \leq \langle D^2
\varphi (x_\varepsilon,t_\varepsilon) p,p \rangle
+\frac{2}{\varepsilon^2} |p-q|^2 +
m\left(\frac{\rho}{\varepsilon^4}\right),
\end{eqnarray}
where $m$ is a modulus of continuity which is independent of
$\rho$ and $\varepsilon.$ In the sequel, $m$ will always
denote a generic modulus of continuity independent of
$\rho$ and $\varepsilon.$

Writing the subsolution viscosity inequality for $U$
and the supersolution inequality for $\tilde{V}$ by
means of the semi-jets and subtracting the inequalities,
we obtain
\begin{eqnarray} \label{eqW}
\varphi_t (x_\varepsilon,t_\varepsilon) &+&
H\big(x_\varepsilon,t_\varepsilon, D\varphi (x_\varepsilon,t_\varepsilon)
+p_\varepsilon , X)
- \mu\, H\big(y_\varepsilon,t_\varepsilon,
\frac{p_\varepsilon}{\mu}, \frac{Y}{\mu}\big)
\nonumber \\
&+& G\big(x_\varepsilon,t_\varepsilon, D\varphi (x_\varepsilon,t_\varepsilon)
+p_\varepsilon , X)
- \mu\, G\big(y_\varepsilon,t_\varepsilon,
\frac{p_\varepsilon}{\mu}, \frac{Y}{\mu}\big)
\leq 0.
\end{eqnarray}

{\it Step 3.}  Estimate of
${\cal{G}}:=
\displaystyle{G\big(x_\varepsilon,t_\varepsilon, D\varphi (x_\varepsilon,t_\varepsilon)
+p_\varepsilon , X)
- \mu\, G\big(y_\varepsilon,t_\varepsilon,
\frac{p_\varepsilon}{\mu}, \frac{Y}{\mu}\big).}$ \\
For sake of simplicity, we set
$$
{c} (x_\varepsilon,t_\varepsilon, \beta)= {c}_x
\ \ \ {\rm and} \ \ \
{c} (y_\varepsilon,t_\varepsilon, \beta)= {c}_y.
$$
We have
\begin{eqnarray*}
{\cal{G}} & = & \mathop{\rm sup}_{\beta\in  {{B}}} \left\{ -\langle
{g}(x_\varepsilon,t_\varepsilon, \beta ), D\varphi
(x_\varepsilon,t_\varepsilon)+ p_\varepsilon \rangle - {f}
(x_\varepsilon,t_\varepsilon, \beta) -{\rm Tr}\left[ {c}_x
{c}_x^T X \right] \right\}
\\ & &
-\mathop{\rm sup}_{\beta\in  {{B}}} \left\{ -\langle
{g}(y_\varepsilon,t_\varepsilon, \beta ), p_\varepsilon \rangle -
\mu {f} (y_\varepsilon,t_\varepsilon, \beta) -{\rm Tr}\left[
{c}_y {c}_y^T Y \right] \right\}
\\ &\geq &
\mathop{\rm inf}_{\beta\in  {{B}}} \left\{ \langle
{g}(y_\varepsilon,t_\varepsilon, \beta )
-{g}(x_\varepsilon,t_\varepsilon, \beta ), p_\varepsilon \rangle -
\langle {g}(x_\varepsilon,t_\varepsilon, \beta ), D\varphi
(x_\varepsilon,t_\varepsilon) \rangle  \right.
\\ & &
\left. - (1-\mu) {f}
(y_\varepsilon,t_\varepsilon, \beta) + {f} (y_\varepsilon,t_\varepsilon, \beta)
- {f} (x_\varepsilon,t_\varepsilon, \beta)
-  {\rm Tr}\left[ {c}_x {c}_x^T X
- {c}_y {c}_y^T Y \right] \right\}
\end{eqnarray*}
From \textbf{(A2)}, if $L_{{g},r}$ is the Lipschitz constant
of ${g}$ in $\overline{B}(\bar{x},r)\times[\bar{t}-r,\bar{t}+r],$
then we have
$$
\langle {g}(y_\varepsilon,t_\varepsilon, \beta )
-{g}(x_\varepsilon,t_\varepsilon, \beta ), p_\varepsilon
\rangle \leq L_{{g},r} |y_\varepsilon -x_\varepsilon| |p_\varepsilon|
\leq 2L_{{g},r} \frac{|y_\varepsilon -x_\varepsilon|^2}{\varepsilon^2}
= m(\varepsilon)\,
$$
and
\begin{eqnarray*}
- \langle {g}(x_\varepsilon,t_\varepsilon, \beta ),
D\varphi (x_\varepsilon,t_\varepsilon) \rangle
\geq -\bar{C}(1+|x_\varepsilon|)
|D\varphi (x_\varepsilon,t_\varepsilon)| \, .
\end{eqnarray*}
By assumption, ${f} \leq 0$ thus
$- (1-\mu) {f} (y_\varepsilon,t_\varepsilon, \beta)\geq 0.$
Again from \textbf{(A2)} it follows that
$$
{f} (y_\varepsilon,t_\varepsilon, \beta)
- {f} (x_\varepsilon,t_\varepsilon, \beta)
\geq -m (|y_\varepsilon -x_\varepsilon|).
$$
Let us denote by $(e_i)_{1\leq i\leq N}$ the canonical basis of $\R^N.$ 
By using (\ref{ineq-fq}), we obtain
\begin{eqnarray*}
{\rm Tr}\left[ {c}_x {c}_x^T X -{c}_y {c}_y^T Y \right]
& = & \sum_{i=1}^N
\langle X{c}_x e_i , {c}_x e_i\rangle
-  \langle Y{c}_y e_i , {c}_y e_i\rangle \\
& \leq &  
{\rm Tr}\left[ {c}_x {c}_x^T D^2 \varphi (x_\varepsilon,t_\varepsilon )
\right]
+\frac{2}{\varepsilon^2} |{c}_x -{c}_y|^2
+m\left(\frac{\rho}{\varepsilon^4}\right)
\\
&\leq &
{\rm Tr}\left[ {c}_x {c}_x^T D^2 \varphi (x_\varepsilon,t_\varepsilon )
\right] +
2L_{{c} ,r}^2 \frac{|x_\varepsilon -y_\varepsilon|^2}{\varepsilon^2}
+ m\left(\frac{\rho}{\varepsilon^4}\right) \\
&\leq &
{\rm Tr}\left[ {c}_x {c}_x^T D^2 \varphi (x_\varepsilon,t_\varepsilon )
\right] +
m(\varepsilon) + m\left(\frac{\rho}{\varepsilon^4}\right),
\end{eqnarray*}
where $L_{{c}, r}$ is a Lipschitz constant for ${c}$ in
$\bar{B}(x,r).$ Hence, since all the modulus are independent of
$\varepsilon,$ $\rho$ and the control, we have
\begin{eqnarray} \label{estimationG}
{\cal{G}} & \geq & -\bar{C}(1+|x_\varepsilon|)
|D\varphi(x_\varepsilon,t_\varepsilon)|+ \mathop{\rm
inf}_{\beta\in  {{B}}} \left\{ -{\rm Tr}\left[
{c}(x_\varepsilon,t_\varepsilon, \beta)
{c}(x_\varepsilon,t_\varepsilon, \beta)^T D^2 \varphi
(x_\varepsilon,t_\varepsilon)\right]
\right\} \nonumber 
+ m(\varepsilon) + m\left(\frac{\rho}{\varepsilon^4}\right).
\end{eqnarray}

{\it Step 4.}  Estimate of
$ \displaystyle{ {\cal{H}} :=
H\big(x_\varepsilon,t_\varepsilon, D\varphi (x_\varepsilon,t_\varepsilon)
+p_\varepsilon , X)
- \mu\, H\big(y_\varepsilon,t_\varepsilon,
\frac{p_\varepsilon}{\mu}, \frac{Y}{\mu}\big).}$ \\
With the same notations as in Step 3, we have
\begin{eqnarray*}
{\cal{H}} &\geq &
\mathop{\rm inf}_{{\alpha}\in \R^k}
\left\{
\langle  {b}(x_\varepsilon,t_\varepsilon, {\alpha} )
- {b}(y_\varepsilon,t_\varepsilon, {\alpha} ), p_\varepsilon \rangle
+
\langle {b}(x_\varepsilon,t_\varepsilon, {\alpha} ),
D\varphi (x_\varepsilon,t_\varepsilon) \rangle  \right. \\
& & \left.+ (1-\mu){\ell}(y_\varepsilon,t_\varepsilon,{\alpha})
+{\ell}(x_\varepsilon,t_\varepsilon,{\alpha})
-{\ell}(y_\varepsilon,t_\varepsilon,{\alpha})
-{\rm Tr}\left[{\sigma}_x {\sigma}_x^T X -{\sigma}_y {\sigma}_y^T Y\right]\right\}.
\end{eqnarray*}
 From \textbf{(A1)}  the following estimates follow:
\begin{eqnarray*}
& & \langle  {b}(x_\varepsilon,t_\varepsilon, {\alpha} )
- {b}(y_\varepsilon,t_\varepsilon, {\alpha} ), p_\varepsilon \rangle
\geq -\bar{C} (1+|{\alpha}|)|x_\varepsilon-y_\varepsilon| |p_\varepsilon|
\geq -\bar{C} |{\alpha}| m(\varepsilon ) + m(\varepsilon )\, ; \\[3mm]
& & \langle {b}(x_\varepsilon,t_\varepsilon, {\alpha} ),
D\varphi (x_\varepsilon,t_\varepsilon) \rangle
\geq
- \bar{C} (1+|x_\varepsilon|)
|D\varphi (x_\varepsilon,t_\varepsilon)|
-  \bar{C}|{\alpha}|
|D\varphi (x_\varepsilon,t_\varepsilon)|\, ; \\[3mm]
& & {\ell}(x_\varepsilon,t_\varepsilon,{\alpha})
-{\ell}(y_\varepsilon,t_\varepsilon,{\alpha}) \geq -(1+|{\alpha}|^2) m_r
(|x_\varepsilon - y_\varepsilon|) \geq -|{\alpha}|^2 m(\varepsilon) +
m(\varepsilon). \\[3mm]
& & (1-\mu){\ell}(y_\varepsilon,t_\varepsilon,{\alpha})\geq
(1-\mu)\left( \frac{\nu}{2}|{\alpha}|^2 +\ell_0 (x, t,{\alpha})\right)
\geq \frac{\nu (1-\mu)}{2}|{\alpha}|^2,
\end{eqnarray*}
where the last inequality follows from the fact that by assumption, $\ell_0(x,t,{\alpha})\geq 0.$
\par
By proceeding exactly as in Step 3 one can show that
\begin{eqnarray*}
{\rm Tr}\left[{\sigma}_x {\sigma}_x^T X -{\sigma}_y {\sigma}_y^T Y\right] \leq
{\rm Tr}\left[{\sigma}_x {\sigma}_x^T
D^2 \varphi (x_\varepsilon,t_\varepsilon)\right]
+m(\varepsilon) + m\left(\frac{\rho}{\varepsilon^4}\right),
\end{eqnarray*}
where $m$ is independent of ${\alpha}.$
Thus, using (\ref{ineg-fondamentale}), we have
\begin{eqnarray}
{\cal{H}} & \geq &
\mathop{\rm inf}_{{\alpha}\in \R^k}\left\{
\left(\frac{\nu (1-\mu)}{2} + m(\varepsilon)\right) |{\alpha}|^2
- \bar{C}(|D\varphi (x_\varepsilon,t_\varepsilon)| + m(\varepsilon))|{\alpha}|
\right\} \nonumber
\\
& & + \mathop{\rm inf}_{{\alpha}\in \R^k}\left\{
-{\rm Tr}\left[ {\sigma}(x_\varepsilon,t_\varepsilon, {\alpha})
{\sigma}(x_\varepsilon,t_\varepsilon, {\alpha})^T
D^2 \varphi (x_\varepsilon,t_\varepsilon)
\right]\right\} \nonumber \\
& &  -\, \bar{C} (1+|x_\varepsilon|)
|D\varphi (x_\varepsilon,t_\varepsilon)|
+ m(\varepsilon ) + m\left(\frac{\rho}{\varepsilon^4}\right)
\nonumber \\
& \geq &
-\frac{(\bar{C} |D\varphi (x_\varepsilon,t_\varepsilon)|
+m (\varepsilon))^2}{2 \nu (1-\mu) + m (\varepsilon)}
- \bar{C} (1+|x_\varepsilon|)
|D\varphi (x_\varepsilon,t_\varepsilon)|
\nonumber \\
&&  \label{estimationH}
+ \mathop{\rm inf}_{{\alpha}\in \R^k}\left\{
-{\rm Tr}\left[ {\sigma}(x_\varepsilon,t_\varepsilon, {\alpha})
{\sigma}(x_\varepsilon,t_\varepsilon, {\alpha})^T
D^2 \varphi (x_\varepsilon,t_\varepsilon)
\right]\right\} +m(\varepsilon)
+ m\left( \frac{\rho}{\varepsilon^4}\right).
\end{eqnarray}

{\it Step 5.} 
Finally, from (\ref{eqW}), (\ref{estimationG}),
(\ref{estimationH}), letting first $\rho$ go to 0 and then
sending $\varepsilon$ to 0,
we obtain
\begin{eqnarray*}
{\cal{L}}[\varphi] (\bar{x},\bar{t}) &= & \varphi_t (\bar{x},\bar{t})
-\frac{\bar{C}^2}{2 \nu (1-\mu)} |D\varphi (\bar{x},\bar{t})|^2
- 2\bar{C} (1+|\bar{x}|)|D\varphi (\bar{x},\bar{t})|
\\
& & +  \!\! \mathop{\rm inf}_{{\alpha}\in \R^k} \!\! \left\{
-{\rm Tr}\left[ {\sigma}(\bar{x},\bar{t}, {\alpha})
{\sigma}(\bar{x},\bar{t}, {\alpha})^T
D^2 \varphi (\bar{x},\bar{t})
\right]\right\} +  \!\!
\mathop{\rm inf}_{\beta\in  {{B}}}\!\! \left\{ -{\rm Tr}\left[
{c}(\bar{x},\bar{t}, \beta) {c}(\bar{x},\bar{t}, \beta)^T
D^2 \varphi (\bar{x},\bar{t})\right] \right\} \leq 0
\end{eqnarray*}
which means exactly that $\Psi$ is a subsolution of (\ref{hjb2}).
~~$\Box$\par

{\it Proof of Lemma \ref{prob-parab}.}
Set $\chi (s,t)= \varphi ({\rm e}^s, t)$ for $(s,t)\in \R\times
[0,+\infty).$
A straightforward calculation shows that
$\varphi$ satisfies (\ref{equa-parab}) if and only if $\chi$ is
a solution of the heat equation
\begin{equation}\label{heat-equation}
\left\{\begin{array}{ll} \chi_t -\chi_{ss}  =0
 & {\rm in}\ \
\R\times (0,T)\\
\chi (s,0)= \varphi_R ({\rm e}^s) & {\rm in}\ \ \R.
\end{array}\right.
\end{equation}
Since the initial data satisfies the growth estimate
$|\chi (s,0)|< {\rm e}^{s^2},$
by classical results on the heat equation (John \cite{john91},
Evans \cite{evans98}),
we know there exists a unique classical solution
$\chi \in C(\R\times [0,T])\times  C^\infty (\R\times (0,T])$
of (\ref{heat-equation}). It is given by the representation
formula: for every $(s,t)\in \R\times [0,T],$
\begin{equation}\label{representaion-form}
\chi (s,t) = \frac{1}{\sqrt{4\pi t}}
\int_{\R} {\rm e}^{-\frac{(s-y)^2}{4t}} \varphi_R ({\rm e}^y)\, dy
= \frac{1}{\sqrt{4\pi t}} \int_{{\rm log}\, R}^{+\infty}
{\rm e}^{-\frac{(s-y)^2}{4t}} ({\rm e}^y -R)\, dy.
\end{equation}
From the above formula, it follows $\chi (s,t)>0$ for all
$(s,t)\in\R\times (0,T].$
Let $h>0.$ We have $\varphi_R ({\rm e}^s)\leq \varphi_R ({\rm e}^{s+h})$
for all $s\in \R.$ Since $\chi (\cdot +h, \cdot)$ is a solution
of (\ref{heat-equation}) with initial data $\varphi_R ({\rm e}^{s+h}),$
by the maximum principle, we obtain
$\chi (s,t)\leq\chi (s+h,t).$ This proves that $\chi$ is nondecreasing
with respect to $s.$
It follows that $\varphi (r,t)=\chi ({\rm log}\, r,t)$ is the
unique solution of (\ref{equa-parab}) and
$\varphi\in C([0,+\infty)\times [0,T])\cap  C^\infty ((0,+\infty)
\times (0,T])$ is positive and nondecreasing.
Moreover, if the initial data is convex, we know that the solution
of a quasilinear equation like (\ref{equa-parab})
is convex in the space variable for every time
(see e.g. Giga {\it et al.} \cite{ggis91}).
 Thus $\varphi(\cdot,t)$
is convex in $[0,+\infty)$ for all $t\in [0,T].$

It remains to prove the estimates (\ref{derivee-phi}).
Noticing that $\varphi_R({\rm e}^s)\leq {\rm e}^s$ and
that $(s,t)\mapsto \varphi_R({\rm e}^s)$ and
$(s,t)\mapsto {\rm e}^{s+t}$ are respectively sub- and supersolution
(in the viscosity sense for example) of (\ref{heat-equation}),
by the maximum principle, we obtain
$\varphi_R({\rm e}^s) \leq \chi (s,t) \leq {\rm e}^{s+t}\leq {\rm e}^{s+T}$
for $(s,t)\in \R\times [0,T].$
It follows
\begin{equation} \label{encadrement1}
\varphi_R(r) \leq \varphi (r,t)\leq {\rm e}^{T} r
\ \ \ {\rm for} \ (r,t)\in [0,+\infty] \times [0,T].
\end{equation}
This gives the first estimate.
To prove the second estimate, we note
that  $\varphi (\cdot,t)$ is a convex nondecreasing
function satisfying (\ref{encadrement1}). It follows
$0\leq \varphi_r (r,t) \leq {\rm e}^{T}$
for $(r,t)\in [0,+\infty] \times (0,T].$
The last assertion is obvious, using the dominated convergence
theorem in (\ref{representaion-form}). It completes the proof of
the lemma.~~$\Box$\par

\section{Applications}
\label{sec-applications}
This Section is divided in two parts. In the first part we consider a
  finite horizon
unbounded stochastic control problem and we characterize the value function
as the unique viscosity solution of the corresponding Dynamic Programming Equation, which is
a particular case of the equation (\ref{hjbintr}). In the second part  we list some
concrete examples of model cases to which the results of Section 1 can be applied.

\subsection{ Unbounded stochastic control problems }

Let  $ (\Omega, {\cal{F}}, ({\cal{F}}_t)_{t\geq 0}, {P})$ be a
filtered probability space, $W_t$ be a ${\cal{F}}_t$-adapted
standard $M$-Brownian motion such that $W_0=0$ a.s. and
let $A$ be a subset of a separable normed space (possibly unbounded).
We consider a {\it finite horizon unbounded stochastic control problem}
for controlled diffusion processes $X_s^{t,x}$ whose dynamic is governed
by a stochastic differential equation of the form
\begin{equation}\label{sode}
\left\{\begin{array}{l}
dX_s^{t,x}=b(X_s^{t,x} ,s,\alpha_s)ds+\sigma(X_s^{t,x} ,s,\alpha_s)dW_s,
\ \ s\in (t,T), \ 0\leq t\leq T,\\
 X_t^{t,x} =x\in \R^N,
\end{array}\right.
\end{equation}
where the control
$\alpha_s\in A,$  $b\colon\R^N\times \R\times A\to\R^N$ is a
continuous vector field, $\sigma$ is a continuous real $N\times M$
matrix.   The pay-off to be minimized is
$$
J(t,x,\alpha)=\mbox{E}_{tx}\{\int_t^T
\ell(X_s^{t,x},s,\alpha_s)\,ds+\psi(X_T^{t,x})\}
$$
where $\mbox{E}_{tx}$ is the expectation with respect to the event $X_t^{t,x}=x$,
the functions $\ell : \R^N\times [0,T]\times A\to \R$ and
$\psi : \R^N\to \R$ are continuous,
$\alpha_s \in{\cal{A}}_t,$ the set of $A$-valued
${\cal{F}}_t$-progressively measurable controls such that
\begin{eqnarray} \label{controleL2}
\mbox{E}_{tx}(\int_t^T|\alpha_s|^2\,ds)<+\infty
\end{eqnarray}
and  $X^{t,x}_s$ is the solution of (\ref{sode}). The value function
is defined by
\begin{equation}
 V(x,t)=\inf_{\alpha_s \in {\cal A}_t} J(t,x,\alpha_s) .
\end{equation}
At least formally,
the Dynamic Programming Equation associated to this control problem is
\begin{equation}\label{stoch}
-\displaystyle\frac{\partial w}{\partial t}
+ \sup_{\alpha\in A}\{-\langle b(x,t,\alpha), D w \rangle -\ell(x,t,\alpha)
-\frac{1}{2} \mbox{Tr}\,(\sigma(x,t,\alpha)\sigma(x,t,\alpha)^T D^2 w )\} =0
\end{equation}
in $\R^N\times(0,T),$
with the terminal value condition $ w (x,T)=\psi(x).$

Our main goal is  to characterize the value function  $V$ as the unique
continuous viscosity solution of $\rec{stoch}$ with the terminal value
condition $V(x,T)=\psi(x).$
We recall that the fact that the value
function is a viscosity solution of the equation  $\rec{stoch}$
is in general obtained  by a direct
use of the {Dynamic Programming Principle}. Since  we are in an
unbounded control framework, the proof of the Dynamic Programming
Principle is rather delicate, thus we  follow another strategy which
consists in comparing directly $V$ with the unique viscosity solution $U$ of
\rec{stoch} obtained by Corollary \ref{thm-existence} when this latter exists.

We make the following assumptions on the data.
\begin{enumerate}
\item[{\bf (S0)}] $A$ is a subset (possibly unbounded)
of a separable complete normed space.

\item[{\bf (S1)}]  $b\in C(\R^N\times [0,T]\times A; \R^N)$ and there
exists $\bar{C}>0$ such that, for all $x,y\in \R^N$,
$t\in[0,T]$ and $\alpha\in A$ we have
\begin{eqnarray*}
|b(x,t,\alpha)-b(y,t,\alpha)|&\le& \bar{C}|x-y|, \\
|b(x,t, \alpha)|&\le& \bar{C}(1+|x|+|\alpha|) ;
\end{eqnarray*}

\item[{\bf (S2)}]  $\sigma\in C(\R^N\times [0,T]\times A; {\cal{M}}_{N,M})$ and there
exists $\bar{C}>0$ such that, for all $x,y\in \R^N$,
$t\in[0,T]$ and $\alpha\in A$ we have
\begin{eqnarray*}
|\sigma(x,t,\alpha)-\sigma(y,t,\alpha)|&\leq& \bar{C}|x-y|, \\
|\sigma(x,t, \alpha)|&\leq& \bar{C}(1+|x|).
\end{eqnarray*}

\end{enumerate}
Moreover we assume that $\ell$ and $\psi$ satisfy respectively
{\bf (A1)}(iii) and {\bf (A3)}.

We first observe that under the current assumptions {\bf (S1)} and
{\bf (S2)} on $b$ and $\sigma,$  for any control
$\alpha\in{\cal{A}}_t$ satifying (\ref{controleL2}) and any random
variable $Z$ such that $\mbox{E} [Z] <\infty,$ there exists a
unique strong solution of the stochastic differential equation \rec{sode}
which satisfies
\begin{eqnarray*}
\mbox{E}\Big\{ \sup_{t\leq s\leq T} |X_s^{t,Z}|^2 \Big\} < \infty
\end{eqnarray*}
(see e.g. Appendix D in \cite{fs93}).
Moreover, we have better estimates on the trajectories of \rec{sode}.
\begin{Lemma}\label{prel2}
Assume {\bf(S0)}, {\bf (S1)} and {\bf(S2)}. For every
$(x,t)\in\R^N\times [0,T]$ and every $ \alpha_s \in {\cal{A}}_t,$
the solution  $X_s^{t,x}$ of \rec{sode} corresponding to
$\alpha_s,$ satisfies the following properties:
\begin{itemize}
\item[(i)] there exists a constant $C>0$ such that
\begin{eqnarray} \label{estimationX2}
{\rm E}_{tx}\Big\{ \sup_{t\leq s\leq T} |X_s^{t,x}|^2 \Big\}
\leq ( |x|^2 + C (T-t) + C{\rm E}_{tx}\int_t^T |\alpha_s|^2\, ds)\, {\rm e}^{C (T-t)}~;
\end{eqnarray}
\item[(ii)] there exists $C_{x,\alpha}>0$ which depends on $x$ and on the
control $\alpha_s$ such that,
for all $s,s^\prime\in [t,T],$
\begin{equation}\label{estraj7}
{\rm E}_{tx}\{ |X_s^{t,x}-X_{s^\prime}^{t,x}| \} \leq
C_{x,\alpha}  |s-s^\prime|^{1/2}.
\end{equation}
In particular for all $\tau \in[t,T]$ we have
\begin{equation}\label{contraj}
{\rm E}_{tx}\Big\{ \sup_{t\leq s\leq\tau} |X_s^{t,x}-x|\Big\} \leq C_{x,\alpha}|\tau|^{1/2}.
\end{equation}
\end{itemize}
\end{Lemma}

{\it Proof of Lemma \ref{prel2}.} We start by proving
(i).  Let us take   an increasing sequence of  $C^2$ functions
$\varphi_R : \R_+ \to \R_+$ such that for all $R>0$, $\varphi^{\prime}_R(r)=0$ if $r>2R$, 
$\varphi^{\prime\prime}_R(r)\le 0$ and $\varphi_R(r)\uparrow r$,
$\varphi^{\prime}_R(r) \uparrow 1$,  as $R\to +\infty.$   By applying Ito's
formula to the process $\varphi_R (|X_s^{t,x}|^2)$, for a.e.
$t\leq \tau\leq T,$ we have (dropping the argument of $\varphi_R$ and its derivatives)
\begin{eqnarray}
&& \hspace*{-0.7cm} \varphi_R (|X_\tau^{t,x}|^2) =  
\varphi_R (|x|^2) + \int_t^\tau
 2\varphi_R^\prime \langle   X_s^{t,x} , b(X_s^{t,x} , s,\alpha_s)\rangle \,ds \nonumber \\
& & \hspace*{-0.7cm} +  \int_t^\tau \!\! \left( \varphi_R^\prime{\rm Tr} [\sigma \sigma^T
(X_s^{t,x} , s, \alpha_s)]+2
\varphi_R^{\prime\prime }|\sigma^T(X_s^{t,x} , s, \alpha_s)X_s^{t,x}|^2 \right)\, ds 
+ \!\! \int_t^\tau  \!\! 2\varphi_R^\prime\langle
 X_s^{t,x} , \sigma (X_s^{t,x} , s, \alpha_s)
dW_s \rangle.
\label{ega-ito1}
\end{eqnarray}
 By using the current assumptions on $b$ and $\sigma$ the following estimate holds
\begin{eqnarray*}\label{estimate}
 & & \int_t^\tau
 \left( 2\varphi_R^\prime \langle   X_s^{t,x} , b(X_s^{t,x} , s,\alpha_s)\rangle
+ \varphi_R^\prime{\rm Tr} [\sigma\sigma^T (X_s^{t,x} , s, \alpha_s)]+2 \varphi_R^{\prime\prime}
|\sigma^T(X_s^{t,x} , s, \alpha_s)X_s^{t,x}|^2 \right)\, ds\\
&& \le 2{C} \int_t^\tau \varphi_R^\prime|X_s^{t,x}| (1+
|X_s^{t,x}|
+ |\alpha_s|)\, ds
 +C\int_t^\tau \varphi_R^\prime (1+ |X_s^{t,x}|^2 ) \, ds \ \
\ {\rm a.s.},
\end{eqnarray*}
where the constant $C$ does not depend neither on the control $\alpha_s$
nor on $R.$ Moreover  we observe that  since $\varphi_R^\prime =0$ for $t>2R$ we have
\begin{eqnarray*}
\mbox{E}_{tx}\Big\{ \int_t^\tau |\varphi_R^\prime \langle
X_s^{t,x}, \sigma (X_s , s, \alpha_s)\rangle|^2 \, ds\Big\}   <
+\infty
\end{eqnarray*}
hence the expectation of the stochastic integral is zero.
By taking the
expectation in (\ref{ega-ito1}) and applying Fubini's Theorem  we
obtain
\begin{eqnarray*}
\mbox{E}_{tx} \{\varphi_R (|X_\tau^{t,x}|^2)\}
\leq  \varphi_R(|x|^2)
+ {C} \int_t^\tau \mbox{E}_{tx}\{\varphi_R^\prime [ 2+ 5 |X_s^{t,x}|^2 + |\alpha_s|^2] \}\, ds.
\end{eqnarray*}
Since $\varphi_R,\varphi_R^\prime  $ are increasing sequences , we can applying Beppo Levi's Theorem.
Therefore by letting $R\to \infty$ we obtain,
 for every $t\leq \tau \leq T,$
\begin{eqnarray*}
\mbox{E}_{tx} \{ |X_\tau^{t,x}|^2\}
& \leq &   
|x|^2+
{C} \int_t^\tau \mbox{E}_{tx}\{ 2+ 5 |X_s^{t,x}|^2 + |\alpha_s|^2 \}\, ds\\
& \le & 
|x|^2+ 5C \int_t^\tau \mbox{E}_{tx} |X_s^{t,x}|^2\,ds
+ C \mbox{E}_{tx}\{ \int_t^\tau |\alpha_s|^2\,ds\} +2C (T-t)
\end{eqnarray*}
Applying Gronwall's Inequality we obtain
\begin{eqnarray*}
\mbox{E}_{tx} \{ |X_\tau^{t,x}|^2 \} \leq
\Big( |x|^2 + 2C(T-t) +  C\, \mbox{E}_{tx} \{ \int_t^T  |\alpha_s|^2\, ds \} \Big)\,
{\rm e}^{5C(T-t)}.
\end{eqnarray*}
We conclude by Doob's maximal inequality, (see, e.g \cite{ks91}).

The proof of (ii) is an extension of the one in the Appendix D in \cite{fs93} and we leave it
to the reader.~~$\Box$

\begin{Proposition} \label{Vlocbornee}
Assume {\bf(S0)}, {\bf (S1)},  {\bf (S2)}, {\bf(A1)}(iii) for $\ell$
and {\bf (A3)} for $\psi.$ Then there exists $0 \leq \tau< T$ such that
the value function $V$ is finite and satisfies the quadratic growth condition
(\ref{croiss-quad}) in $\R^N\times [\tau, T]$.
\end{Proposition}

{\it Proof of Proposition \ref{Vlocbornee}.}
We aim at showing that that if the constants $\rho, C>0$ are large enough then there exists
$\tau>0$ depending on $\rho$ such that $|V(x,t)|\le C(1+|x|^2){\rm e}^{\rho(T-t)}$
for all $(x,t)\in \R^N\times [\tau,T].$
The upper estimate is obtained by majorizing directly the value
function with the cost functional corresponding to a constant
control and using  the estimates of the trajectories  in  Lemma
\ref{prel2}.
The most difficult is to prove the estimate from below since $V$ is defined
by an infimum.
To this purpose we take any control $\alpha_s\in {\cal{A}}_t$. By applying Ito's formula to
the process $(1+|X_s^{t,x}|^2){\rm e}^{\rho(t-s)}$, $X_s^{t,x}$ being the trajectory corresponding
to $\alpha_s$, we have the following estimate
\begin{eqnarray}\nonumber
 d[(1+|X_s^{t,x}|^2){\rm e}^{\rho(T-s)}]
& = & -\rho {\rm e}^{\rho(T-s)}(1+|X_s^{t,x}|^2)ds
+ {\rm e}^{\rho(T-s)}\mbox{Tr}\,
( \sigma \sigma^T (X_s^{t,x},s,\alpha_s))ds \\
& & +  2 {\rm e}^{\rho(T-s)}\langle X_s^{t,x}, b(X_s^{t,x},s,\alpha_s)ds
+ \sigma(X_s^{t,x},s,\alpha_s) dW_s\rangle.
\label{estraj1}
\end{eqnarray}
Integrating both sides of \rec{estraj1} from $t$ to $T$ and taking the expectation we get:
\begin{eqnarray} \label{estraj2}&  & \mbox{E}_{tx}\{ 1+|X_T^{t,x}|^2 \} - (1+|x|^2) {\rm e}^{\rho(T-t)}
\\ & = &
\mbox{E}_{tx} \Big\{ \int_t^T   \Big( - \rho (1+|X_s^{t,x}|^2)
+ 2 \langle X_s^{t,x} ,  b(X_s^{t,x},s,\alpha_s) \rangle
+ \mbox{Tr} (\sigma \sigma^T (X_s,s,\alpha_s)\Big)  {\rm e}^{\rho(T-s)} \, ds  \Big\}.\nonumber
\end{eqnarray}
We notice that in the above estimate we  supposed   that the
expectation of the stochastic integral is zero. This is false in
general but we can overcome such a difficulty by an
approximation argument  which is similar to the one used in the proof of Lemma \ref{prel2}.
Now for any $\varepsilon$-optimal control $\alpha_s$ for $V(x,t),$ by using \rec{estraj2}, we get
\begin{eqnarray*}
V(x,t)+C(1+|x|^2){\rm e}^{\rho(T-t)}+\varepsilon 
 & \ge &  \mbox{E}_{tx}
\Big\{ \int_t^T \Big( \ell(X_s^{t,x},s,\alpha_s)
 -2C {\rm e}^{\rho(T-s)}
\langle X_s^{t,x}, b(X_s^{t,x},s,\alpha_s)\rangle\\& &  \quad
-C{\rm e}^{\rho(T-s)}\mbox{Tr} (\sigma \sigma^T (X_s^{t,x},s,\alpha_s))
  +\rho {\rm e}^{\rho(T-s)}C(1+|X_s|^2) \Big) \,ds\\ & &\quad
+\psi(X_T^{t,x})+C(1+|X_T^{t,x}|^2)  \Big\} \\
 &=& \mbox{E}_{tx}\{\int_t^T
\bar \ell(X_s^{t,x},s,\alpha_s)\,ds+\bar\psi(X_T^{t,x}) \},
\end{eqnarray*}
where
$\bar\ell(x,t,\alpha):= 
\ell(x,t, \alpha)-2C{\rm e}^{\rho (T-t)}\langle
b(x,t,\alpha), x\rangle -2C{\rm e}^{\rho (T-t)}\mbox{Tr}\,a(x,t,\alpha) 
+ C\rho {\rm e}^{\rho (T-t)}(1+|x|^2)$
and $\bar\psi(x):=\psi(x)+C(1+|x|^2).$
By analogous arguments as those used in Section \ref{sec-1-comp}
one can see that for $\rho, C>0$ large enough there is
$\tau>0$ such that $\bar\ell$ and $\bar\psi$ are nonnegative in $\R^N\times[\tau,T].$
Thus we can conclude since $\varepsilon $ is
arbitrary. ~~$\Box$

\begin{Remark}\rm
If $\ell,\psi$ are bounded from below, namely they satisfy, for some $C>0,$ the following two
conditions
$\ell(x,t,\alpha) \ge \nu|\alpha|^2-C$ and $\psi(x) \ge -C,$
then $V$ is finite and satisfies the growth condition (\ref{croiss-quad}) in $\R^N\times [0,T]$
(i.e. for all time).
\end{Remark}

Next we prove that $V$ is the unique viscosity solution of \rec{stoch}.
We start with the following Proposition.
\begin{Proposition}\label{supersolution}
Under the assumptions of Proposition \ref{Vlocbornee} we have \\
(i) (Super-optimality principle) for all $t\in (\tau,T]$ and $0<h\le T-t$ and for all stopping time
$t\leq \theta \leq T,$ we have
\begin{equation}
 V(x,t)\ge \inf_{\alpha_s \in {\cal A}_t}{\rm E}_{tx}\{\int_t^{(t+h)\land\theta}
 \ell(X_s^{t,x},s,\alpha_s)\,ds+V_* ( X_{(t+h)\land\theta}^{t,x},(t+h)\land\theta)\}.
\end{equation}
(ii) the function $V$ is a supersolution of \rec{stoch} in $\R^N\times [\tau,T].$
\end{Proposition}

{\it Proof of Proposition \ref{supersolution}.} The
proof of (i)  is a standard routine  and  we refer the reader for
instance  to the book of Yong and Zhou \cite{yz99}. The opposite inequality is more
delicate, see Krylov \cite{krylov80, krylov01}.

We turn to the proof of (ii),
showing that the super-optimality principle  implies that $V_*$ is
a viscosity supersolution of \rec{stoch}. Let $\phi\in
C^{2}(\R^N\times [0,T])$ and $(\bar x,\bar t) \in \R^N\times (\tau
, T)$ be a local minimum of $V_* -\phi.$ We can assume that
$V_*(\bar x,\bar t)=\phi(\bar x,\bar t)$ and that the maximum is
strict, i.e. $V_*(x,t)> \phi(x,t)$ for all $(x,t)\in \bar B(\bar x,\varepsilon)\times
[\bar t -\varepsilon , \bar t +\varepsilon ]$  with
$(x,t)\not=(\bar x,\bar t)$ (see \cite{bcd97} or \cite{barles94}   ).
We assume by contradiction that there exists $\delta_\varepsilon >0$
such that for all $(x,t)\in \bar B(\bar x,\varepsilon)\times [\bar
t -\varepsilon , \bar t +\varepsilon ],$ we have
\begin{equation}\label{stoch2}
-\phi_t(x,t)+\sup_{\alpha\in A} \{ -\langle b(x,t,\alpha), D\phi (x,t)\rangle -\ell(x,t,\alpha)
-{\rm Tr}\,[\frac{1}{2}\sigma \sigma^T(x,t,\alpha)D^2 \phi (x,t) ] \}
\leq -\delta_\varepsilon.
\end{equation}
Since $(\bar x,\bar t)$ is a strict minimum of $V_* -\phi,$ it follows that there exists
$\eta_\varepsilon$ such that
\begin{equation}\label{bord-boule}
V_*(x,t)\geq \phi(x,t)+\eta_\varepsilon \ \ \
{\rm for \ all} \ (x,t)\in \partial B(\bar x,\varepsilon)\times
[\bar t -\varepsilon , \bar t +\varepsilon ].
\end{equation}
From now on, we fix $0<h<\varepsilon /2$ such that $h\delta_\varepsilon < \eta_\varepsilon.$
Let us denote by $\tau_{t,x}$ the exit time of the trajectory $X^{t,x}_s$ from the ball
$B(\bar x,\varepsilon).$ We first observe that by the continuity of the trajectory
(see Lemma \ref{prel2}), we have $\tau_{t,x}>t$ for all $(x,t)\in B(\bar x,\varepsilon)\times
[0,T).$
For every $(x,t)\in B(\bar x,\varepsilon)\times (\bar t -\varepsilon /2, \bar t +\varepsilon /2),$
there exists a control $\alpha_s \in {\cal{A}}_t$ such that
\begin{eqnarray*}
V(x,t) +\frac{\delta_\varepsilon h}{2} \geq \mbox{E}_{tx}
\Big\{
\int_t^{(t+h)\land\tau_{t,x}} \ell(X^{t,x}_s,s,\alpha_s)\,ds
+V_* ( X^{t,x}_{(t+h)\land\tau_{t,x}},(t+h)\land\tau_{t,x})
\Big\}.
\end{eqnarray*}
Since $V_*\geq \phi$ in
$\bar B(\bar x,\varepsilon)\times  [\bar t -\varepsilon , \bar t +\varepsilon ],$ if
$\tau_{t,x} < t+h,$ then, from (\ref{bord-boule}), we have
$$
V_* ( X^{t,x}_{(t+h)\land\tau_{t,x}},(t+h)\land\tau_{t,x}) \geq \phi (X^{t,x}_{\tau_{t,x}}, \tau_{t,x})
+\eta_\varepsilon.
$$
Therefore the following estimate holds
\begin{eqnarray} \nonumber
V(x,t) +\frac{\delta_\varepsilon h}{2} & \geq & \mbox{E}_{tx}
\Big\{ \Big[
\int_t^{\tau_{t,x}} \ell(X^{t,x}_s,s,\alpha_s)\,ds
+\phi (X^{t,x}_{\tau_{t,x}}, \tau_{t,x}) + \eta_\varepsilon \Big]
\11_{\{\tau_{t,x}<t+h\}}
\Big\} \\
\nonumber
& & \hspace*{0.5cm} +
\mbox{E}_{tx} \Big\{ \Big[
\int_t^{t+h} \ell(X^{t,x}_s,s,\alpha_s)\,ds
+\phi (X^{t,x}_{t+h}, t+h) \Big]
\11_{\{\tau_{t,x}\geq t+h\}}\Big\} \\
\label{longest}
&\geq &
\mbox{E}_{tx} \Big\{ [I(\tau_{t,x}) + \eta_\varepsilon] \11_{\{\tau_{t,x}<t+h\}} \Big\}
+ \mbox{E}_{tx} \Big\{ I(t+h) \11_{\{\tau_{t,x}\geq t+h\}} \Big\},
\end{eqnarray}
where for all $\tau^\prime>0$
\begin{eqnarray*}
I(\tau^\prime)=
\int_t^{\tau^\prime} \ell(X^{t,x}_s,s,\alpha_s)\,ds + \phi (X^{t,x}_{\tau^\prime}, \tau^\prime).
\end{eqnarray*}
Applying Ito's formula to the process $\phi (X^{t,x}_{\tau^\prime}, \tau^\prime),$ we obtain
\begin{eqnarray*}
I(\tau^\prime)&=&
\int_t^{\tau^\prime} \Big(
\ell(X^{t,x}_s,s,\alpha_s) + \phi_t (X^{t,x}_s,s) +
\langle D\phi (X^{t,x}_s,s), b(X^{t,x}_s,s,\alpha_s) \rangle\\
&&+\frac{1}{2} {\rm Tr}\, [ \sigma \sigma^T (X^{t,x}_s,s,\alpha_s)D^2\phi]
\Big)\, ds 
+ \phi (x,t)
+ \int_t^{\tau^\prime} \langle D\phi (X^{t,x}_s,s), \sigma (X^{t,x}_s,s,\alpha_s) dW_s \rangle
\ \ \ {\rm a.s.}
\end{eqnarray*}
Note that the expectation of the above stochastic integral is zero for $\tau^\prime \in [t,t+\varepsilon].$

Now we can estimate the two last terms in (\ref{longest}).
For the first term, we have
\begin{eqnarray*}
& & \mbox{E}_{tx} \Big\{ [I(\tau_{t,x}) + \eta_\varepsilon] \11_{\{\tau_{t,x}<t+h\}} \Big\}
\\ & \geq &
-\mbox{E}_{tx}\Big\{ \Big[
\int_t^{\tau_{t,x}}\Big( -\phi_t(X^{t,x}_s,s)+ \sup_{\alpha\in A}
\{ - \ell(X^{t,x}_s,s,\alpha ) -
\langle D\phi (X^{t,x}_s,s), b(X^{t,x}_s,s,\alpha ) \rangle
\\ & & \hspace*{3.5cm}
- \frac{1}{2} {\rm Tr}\, [ \sigma \sigma^T (X^{t,x}_s,s,\alpha )D^2\phi] \}
\Big)\, ds
-\phi(x,t) - \eta_\varepsilon \Big] \11_{\{\tau_{t,x}<t+h\}}
\Big\}.
\end{eqnarray*}
Since $X^{t,x}_s \in B(\bar x,\varepsilon)$ when $s\leq \tau_{t,x}$
and since $t+h < \bar t + \varepsilon,$ from (\ref{stoch2}), we get
\begin{eqnarray} \nonumber
\mbox{E}_{tx} \Big\{ [I(\tau_{t,x}) + \eta_\varepsilon] \11_{\{\tau_{t,x}<t+h\}} \Big\}
& \geq &
- \mbox{E}_{tx} \Big\{ \Big[ \int_t^{\tau_{t,x}} (-\delta_\varepsilon)\, ds
-\phi(x,t) - \eta_\varepsilon  \Big] \11_{\{\tau_{t,x}<t+h\}} \Big\}
\\ \nonumber
& \geq &  \delta_\varepsilon  \mbox{E}_{tx}\left[(\tau_{t,x} -t)\11_{\{\tau_{t,x}<t+h\}}\right]
+ (\eta_\varepsilon + \phi(x,t) ) P(\{\tau_{t,x}<t+h\} ) \nonumber \\
&\geq & (\eta_\varepsilon + \phi(x,t) ) P(\{\tau_{t,x}<t+h\} ).  \label{estima1}
\end{eqnarray}
For the second term, we proceed in the same way, noting that, if $\tau_{t,x}\geq t+h,$
then for all $t\leq s \leq t+h,$ $X^{t,x}_s \in B(\bar x,\varepsilon)$ and
it allows us to apply (\ref{stoch2}). More precisely we have
\begin{eqnarray} 
& & \mbox{E}_{tx} \Big\{I(t+h)\11_{\{\tau_{t,x}\geq t+h\} } \Big\} \label{estima2}\\
&\geq &
-\mbox{E}_{tx}\Big\{ \Big[
\int_t^{t+h} \Big( -\phi_t(X^{t,x}_s,s)+ \sup_{\alpha\in A}
\{ - \ell(X^{t,x}_s,s,\alpha ) -
\langle D\phi (X^{t,x}_s,s), b(X^{t,x}_s,s,\alpha ) \rangle
\nonumber\\  & & \hspace*{4.5cm}
- \frac{1}{2} {\rm Tr}\, [ \sigma \sigma^T (X^{t,x}_s,s,\alpha )D^2\phi] \}
\Big)\, ds
-\phi(x,t) \Big] \11_{\{\tau_{t,x}\geq t+h\}}
\Big\} \nonumber
 \\
& \geq &
( \delta_\varepsilon  h + \phi(x,t) ) P(\{\tau_{t,x}\geq t+h\}). \nonumber
\end{eqnarray}
Combining (\ref{longest}), (\ref{estima1}) and  (\ref{estima2}), we get
\begin{eqnarray*}
V(x,t) +\frac{\delta_\varepsilon h}{2} & \geq &
  \eta_\varepsilon  P(\{\tau_{t,x}<t+h\} )
+ \delta_\varepsilon  h \, P(\{\tau_{t,x}\geq t+h\}) \\
&& 
+  \phi(x,t) \left[P(\{\tau_{t,x}<t+h\} ) + P(\{\tau_{t,x}\geq t+h\})\right].
\end{eqnarray*}
Since $\eta_\varepsilon > \delta_\varepsilon h$
and $P(\{\tau_{t,x}<t+h\} ) + P(\{\tau_{t,x}\geq t+h\})=1,$
we get
$$
V(x,t)\geq \phi( x,t) + \frac{\delta_\varepsilon h}{2}.
$$
The above inequality is valid for all $(x,t)\in
B(\bar x,\varepsilon)\times (\bar t -\varepsilon /2, \bar t +\varepsilon /2),$
thus we have
$$
\liminf_{(x,t)\to(\bar x,\bar t)}V(x,t) = V_*(\bar x,\bar t)
\geq \phi(\bar x,\bar t)+{\frac{\delta_\varepsilon h}{2}}
$$
which is a contradiction with the choice of $\phi.$ ~~$\Box$

\begin{Theorem} \label{valeur-unique-sol}
Under the assumptions of Proposition \ref{Vlocbornee},
the function $V$ is the unique continuous viscosity solution of $\rec{stoch}$
in $\R^N\times [\tau ,T].$
\end{Theorem}

{\it Proof of Theorem \ref{valeur-unique-sol}.}
Let $U$ be the unique solution of $\rec{stoch}$ in $\R^N\times [\tau,T]$ such that
$U(x,T)=\psi(x)$ given by Theorem \ref{thm-existence}. Our goal is to  prove that $V\equiv U$
in $\R^N\times[\tau ,T].$
The inequality $U\le V_*$  follows by combining
Proposition \ref{supersolution} and Theorem \ref{thm-unicite}.
To show that $V^*\le U$ in $\R^N\times [\tau ,T],$  we proceed as
follows:

{\it Step 1.} We consider the functions $\widetilde
 V(x,t):=V(x,t)-C(1+|x|^2){\rm e}^{\rho (T-t)}$ and $\widetilde
 U(x,t):=U(x,t)-C(1+|x|^2){\rm e}^{\rho(T- t)}.$ As it is proved in Step 2
of the proof of Theorem \ref{thm-unicite}, $\widetilde U$ is the unique
solution of
\begin{eqnarray}
\left\{\begin{array}{l}
\displaystyle
-w_t+\sup_{\alpha\in A}\{-\frac{1}{2}\mbox{Tr}(\sigma \sigma^T(x,t,\alpha)D^2
w)-\langle b(x,t,\alpha),  Dw\rangle -\bar\ell(x,t,\alpha)\}=0 
\ \ {\rm in} \  \R^N\times (\tau,T), \\
w (x,T)=\bar\psi (x),
\end{array}\right. \label{equa-preuve}
\end{eqnarray}
where
$\bar\ell(x,t,\alpha) 
:= \ell(x,t, \alpha)+2C{\rm e}^{\rho (T-t)}\langle
b(x,t,\alpha), x\rangle +C{\rm e}^{\rho (T-t)}\mbox{Tr}\,\sigma \sigma^T(x,t,\alpha)
-C\rho {\rm e}^{\rho (T-t)}(1+|x|^2)$
and $\bar\psi(x):=\psi(x)-C (1+|x|^2).$

{\it Step 2.}
Claim:  for all $(x,t) \in \R^N\times [\tau, T],$ $\widetilde V$ satisfies
 \begin{equation}\label{esttildev}
 \widetilde
 V(x,t)\leq \inf_{\alpha_s \in {\cal A}_t} \mbox{E}_{tx}\{\int_t^T
\bar \ell(X_s,s,\alpha_s)\,ds+\bar\psi(X_t) \}.
\end{equation}
To prove the claim let us take  any $\alpha_s \in {\cal{A}}_t.$
Arguing exactly as in the proof of Proposition \ref{Vlocbornee}, from (\ref{estraj2}), we have  
\begin{eqnarray*}
\widetilde
 V(x,t)&\leq& \mbox{E}_{tx}\Big\{
\int_t^T \Big(
 \ell(X^{x,t}_s,s,\alpha_s)+ 2 {\rm e}^{\rho(T-s)} \langle X_s^{x,t}, b(X^{x,t}_s,t,\alpha_s)\rangle \\
&& \hspace*{0.6cm}+ {\rm e}^{\rho(T-s)}\, \mbox{Tr} (\sigma \sigma^T(X^{x,t}_s,s,\alpha_s))
-\rho {\rm e}^{\rho(T-s)}(1+|X^{x,t}_s|^2)\Big) \,ds
 +\psi(X^{x,t}_T)-C(1+|X^{x,t}_T|^2) \Big\}\\
&=& \mbox{E}_{tx}\{\int_t^T
\bar \ell(X^{x,t}_s,s,\alpha_s)\,ds+\bar\psi(X^{x,t}_T) \}.
\end{eqnarray*}
Since $\alpha_s$ is arbitrary we get \rec{esttildev} and prove the claim.

{\it Step 3.}
Choose $C, \rho>0$ large enough  so that, for all $(x,t,\alpha)\in\R^N\times[\tau,T]\times A,$ we have
\begin{eqnarray*}
 -\widetilde{C} {\rm e}^{2\rho (T-t)}(1+|x|^2)+
\frac{\nu}{4}|\alpha|^2\le \bar\ell(x,t,\alpha)&\le& -{C}{\rm e}^{\rho (T-t)}(1+|x|^2)+C {\rm e}^{\rho
(T- t)}(1+|\alpha|^2),\\
  -2C(1+|x|^2)\leq \bar\psi(x)&\le& 0,
\end{eqnarray*}
where $\widetilde{C}$ depends only on $C$ and $\nu.$
For all real $R>0$ and all integer $n>0,$ we set
$A_n:=\{\alpha\in A : |\alpha|\le n\},$
$\bar\ell_R(x,t,\alpha):= {\rm max}\{\bar\ell(x,t,\alpha), -R\}$ and
$\bar\psi_R(x,t):={\rm max}\{ \bar\psi(x,t), -R \}.$
We observe
that $\bar \ell_R\colon\R^N\times[\tau ,T]\times A_n\to\R$  is  bounded and uniformly
continuous in $x\in\R^N$ uniformly with respect to $(t,\alpha)\in [\tau,T]\times A_n$
and $\psi_R\colon\R^N\to\R $ is bounded and uniformly continuous in $\R^N.$
Set
\begin{eqnarray*}
\bar{H}(x,t,p,X):=\sup_{\alpha\in A}
\{-\frac{1}{2}\mbox{Tr}(\sigma\sigma^T(x,t,\alpha)X)-\langle b(x,t,\alpha), p\rangle-\bar\ell(x,t,\alpha)\}, \\
H^R_n(x,t,p,X):=\sup_{\alpha\in A_n}
\{-\frac{1}{2}\mbox{Tr}(\sigma\sigma^T(x,t,\alpha)X)-\langle b(x,t,\alpha), p\rangle-\bar\ell_R(x,t,\alpha)\}
\end{eqnarray*}
and define
$$
V^R_n(x,t)=\inf_{\alpha_s\in {\cal A}_t^n}
\mbox{E}_{tx}\left\{\int_t^T \bar
\ell_R(X^{x,t}_s,s,\alpha_s)\,ds+\bar\psi_R(X_T^{x,t}) \right\},
$$
where ${\cal{A}}_t^n$ is the set of $A_n$-valued
${\cal{F}}_t$-progressively measurable controls such that
(\ref{controleL2}) holds.
The function $V^R_n$ is now the value function of a stochastic control problem with
bounded controls and uniformly continuous datas $b,$ $\sigma,$ $\bar\ell_R$ and $\bar\psi_R.$
These assumptions enter the framework of Yong and Zhou \cite{yz99}. We deduce that $V^R_n$
is the unique continuous viscosity solution of
\begin{equation}
-\frac{\partial V^R_{n}}{\partial t}+ H^R_n(x,t,D V^R_n,D^2V^R_n)=0 \ \ {\rm in} \ \R^N\times(\tau,T),
\end{equation}
with terminal condition $V^R_{n}(x,T)=\bar\psi_R (x)$ in $\R^N.$
Moreover for all compact subsets $K\subset \R^N$  there exists $M>0$
independent on $R$ and $n$  such that
$||V^R_n||_\infty\leq M$ in $K\times [\tau ,T].$
Indeed, take any constant control $\alpha_s=\bar\alpha \in A_1,$ by definition of $V^R_n$ for all
$R, n>0$ and for every $(x,t)\in\R^N\times[\tau,T]$  we have
\begin{eqnarray*}
V^R_n(x,t)\le \mbox{E}_{tx}\{\int_t^T \bar
\ell_R(X^{x,t}_s,s,\bar{\alpha})\,ds+\bar\psi_R(X_T^{x,t}) \}
\leq  \int_t^T C {\rm e}^{\rho(T-s)}(1+\bar\alpha^2)\,ds
\leq \frac{C}{\rho} {\rm e}^{\rho(T-t)}(1+\bar\alpha^2),
\end{eqnarray*}
and on the other hand we have $V^R_n(x,t)\ge \widetilde V(x,t)\geq
-\widetilde C{\rm e}^{2\rho(T-t)}(1+|x|^2).$

Finally  one can readily see that  $H_n^R$ converges locally uniformly to $\bar{H}$ as
$n,R\to\infty.$
Thus by
applying the half-relaxed limits method (see  Barles and Perthame \cite{bp88}), the functions
$$
\overline V(x,t) = {\limsup}^* V_n^R(x,t)
=\limsup_{{\displaystyle{\mathop{\scriptstyle{(y,s)\to (x,t)}}_{ n,R\to +\infty}}}}V_n^R(y,s)
\ \ \ {\rm and} \ \ \ \underline{V}(x,t)=
{\liminf}^* V_n^R(x,t)=\liminf_{{\displaystyle{\mathop{\scriptstyle{(y,s)\to
(x,t)}}_{n, R\to +\infty}}}}V_n^R(y,s)
$$
are respectively viscosity sub and supersolution of \rec{equa-preuve}.
Theorem \ref{thm-unicite} yields
$\overline V(x,t)\le \widetilde U(x,t)\le \underline V(x,t).$
On the other hand by construction we have also
$\widetilde{V}^*(x,t)\le{\limsup}^*V^R_n(x,t).$
It follows $\widetilde{V}^*(x,t)\le \widetilde U(x,t)$
and we can conclude. ~~$\Box$\par

\subsection{Some examples}

\begin{Example} \rm \label{stochastic-LQ}
A model case we have in mind is the so-called stochastic linear regulator problem which is a
stochastic perturbation of the deterministic linear quadratic problem. In this case,
the stochastic differential (\ref{sode}) is linear and reads
$$
dX_s^{t,x}=[B(s)X_s^{t,x}+C(s)\alpha_s]ds+\sum_j [ C_j(s)X_s^{t,x} +D_j(s)]  dW^j_s
$$
and the expected total cost to minimized is
$$
J(x,t,\alpha_s )=\mbox{E}_{tx}\{\int_t^T
[\langle X_s^{t,x},  Q(s)X_s^{t,x}\rangle +\langle \alpha_s, R(s)\alpha_s\rangle  ]\,ds
+\langle X_T^{t,x},  G X_T^{t,x} \rangle \}.
$$
The previous results apply if
the functions $B(\cdot),C(\cdot), C_j(\cdot), D_j(\cdot), Q(\cdot), R(\cdot)$ and $G$ are
deterministic continuous matrix-valued functions of suitable size
and if $R(s)$ is a positive definite symmetric matrix.
Deterministic and stochastic linear quadratic problems were extensively studied.
For a survey we refer for instance to the books of Bensoussan \cite{bensoussan82},
Fleming and Rishel \cite{fr75}, Fleming and Soner \cite{fs93}, {\O}ksendal \cite{oksendal98},
Yong and Zhou \cite{yz99} and references therein.
\end{Example}

\begin{Example} \rm \label{exple-math-finance}
Equations of the type  \rec{stoch} are largely considered in
mathematical finance. See the introductory books quoted in the introduction
or Pham \cite{pham02}. In particular recently Benth and Karlsen \cite{bk03} studied the following
semilinear elliptic partial equation:
\begin{equation}\label{stoch1}
- w_t-\frac{1}{2}\beta^2  w_{xx}+F(x, w_x)=0~~\quad\mbox{in $\R\times(0,T),$}
\end{equation}
with the final condition
$ w(x,T)=0.$
The nonlinear function $F$ is given by
$$
F(x,p)=\frac{1}{2}\delta^2 p^2-\{\alpha(x)
-\frac{\mu(x)\beta\rho}{\sigma(x)}\} p-\frac{\mu^2(x)}{\sigma^2(x)}
$$
where $\delta, \beta, \rho$ are real constant and
$\alpha(x)$ and $\displaystyle\frac{\mu}{\sigma}(x)$ are $C^1$ functions satisfying
$$|\alpha(x)|,| \frac{\mu}{\sigma}(x)|\le C|x|,\quad\mbox{for all $x\in\R.$}
$$
Their main motivation is to determine via the solution of
(\ref{stoch1})   the minimal entropy martingale measure in
stochastic market. The conditions they assume on data fall within
the assumptions of Section \ref{sec-1-comp}.
\end{Example}

\begin{Example} \rm \label{example-risk-sensitive}
Another application of the results obtained in Section \ref{sec-1-comp}
is given by the finite time-horizon \rsen limit problem for nonlinear systems.
In the stochastic risk-sensitive problem, (\ref{sode}) reads
\begin{eqnarray}\label{rsdyn}
 dX^{t,x, \varepsilon}_s=g(X^{t,x, \varepsilon}_s,\beta_s)\,dt+ {\sqrt{\frac{\eps}{\gamma^2}}}
c(X^{t,x,\varepsilon}_s)\,dW_s,
\end{eqnarray}
where $X^{t,x, \varepsilon}_s\in \R^N$ depends on the parameter $\varepsilon >0,$
$g$ represents the nominal dynamics with control $\beta_s\in B,$ a compact normed space
and $c$ is an $N\times k$-valued diffusion coefficient, $\varepsilon$ is a measure of
the risk-sensitivity and $\gamma$ is the disturbance attenuation level.
The cost criterion is of the form
\begin{equation}
J^\eps(x,t,\beta_s)=\expect\exp\left\{
\frac{1}{\eps}\left[\int_t^T
f(X_s^{t,x,\varepsilon},\beta_s)\,dt+\psi(X_T^{t,x,\varepsilon})\right]\right\}
\end{equation}
and the value function is
\begin{equation}
\ba{rl}
\vvep(x,t)&=\inf\limits_{\beta_s\in{\cal{B}}_t}\eps\log\jjep(x,t,\beta_s)
 =\eps\log\inf\limits_{\beta_s\in{\cal{B}}_t }\jjep(x,t,\beta_s), \ea
\end{equation}
where ${\cal{B}}_t$ is the set of $B$-valued, \ftpm controls such that there exists a strong
solution to \er{rsdyn}. The dynamic programming equation associated to this problem is
\begin{equation}\label{rspde}
\left\{\begin{array}{ll}-\displaystyle{\frac{\partial  w}{\partial t}}
-\displaystyle\invtwogam \langle D  w, a(x)D  w \rangle +\tilde{G}(x,D w)
-\displaystyle\eptwogam\mbox{Tr}(a(x)D^2 w )
 =0
&\mbox{in $\R^N\times[\tau,T],$}  \\[3mm]
 w(x,T)=\psi(x), &
\end{array}\right.
\end{equation}
where $a(x) = c(x)c(x)^T$ and
$$
\tilde{G}(x,p)=\max_{\beta\in B}
\{\langle -g(x,\beta), p\rangle-f(x,\beta)\}.
$$
We note that
\begin{equation}
- \invtwogam \langle p, a(x)p\rangle =\min_{\alpha\in\R^k}
\left\{ \gamovertwo |\alpha|^2 - \langle c(x) \alpha, p\rangle \right\}.
\end{equation}
In \cite{dlme02} it is shown that as  $\eps\downarrow 0$
(i.e. as the problem becomes infinitely risk averse),  the value
function of the risk-sensitive problem converges to that of a $H_\infty$
robust control problem. This problem can be considered
as a differential game with the following cost functional:
\begin{equation}
J(x,t,\alpha ,\beta)=\int_t^T
\Big( f(y_x(s),\beta_s)-{\frac{\gamma^2}{ 2}}|\alpha_s|^2\Big)
ds+\psi(y_x(T)),
\end{equation}
where $\alpha_s \in {\cal{A}}_t:= L^2([t,T],\R^k)$ is the control
of the maximizing player,
$\beta_s\in {\cal{B}}_t :=\{ {\rm measurable \ functions} \ [t,T]\to B \}$ is the
control of the minimizing player, and
$y_x(\cdot)$ is the unique solution of the following dynamical
system:
$$
\left\{
\begin{array}{l}
y^\prime (s) = g(y(s),\beta_s)+c (y(s)) \alpha_s, \\
y(t)=x.
\end{array}
\right.
$$
Note that we switch notation from $X_t$ to $y(t)$ to emphasize that
the paths are now (deterministic) solutions of ordinary differential equations
rather than (stochastic) solutions of stochastic differential equations.
The dynamic programing equation associated to the robust control problem is a
first order equation given by
\begin{equation}\label{robust}
-\displaystyle\frac{\partial  w}{\partial t}+\min_{\alpha\in\R^k}
\left\{\frac{\gamma^2}{2}|\alpha|^2-\langle c(x)\alpha, D w\rangle\right\}
+\tilde G(x,D w)=0~~\mbox{in $ \R^N\times (\tau,T)$}
\end{equation}
with the terminal condition
$ w(T,x)=\psi(x).$
One of the key tools to get this convergence result is the uniqueness
property for (\ref{robust}). In Da Lio and McEneaney \cite{dlme02},
the authors characterized the value function of the $H_\infty$ control as the unique solution
to (\ref{robust}) in the set of locally Lipschitz continuous functions
growing at most quadratically with respect to the state variable.
Moreover the uniqueness result in \cite{dlme02} is obtained by using representation formulas
of locally Lipschitz solutions of \rec{robust}.
We remark that the Comparison Theorem \ref{thm-unicite} not only
improves the uniqueness result for (\ref{robust})
obtained in  \cite{dlme02} (in the sense it holds in a larger class of functions)  but it also
should allow us to  prove  the convergence result in \cite{dlme02} under weaker assumptions, only
the equi-boundedness estimates of the solutions of \rec{rspde} being enough  by means of
the half-relaxed limit method.
\end{Example}

\section{Study of related equations}
\label{sec-related}

In this Section we focus our attention to
Hamilton-Jacobi equations of the form
\begin{equation}\label{hjb-related} \left\{\begin{array}{ll}
\displaystyle\frac{\partial  w}{\partial t}
+ \langle \Sigma (x,t) D w,D w\rangle
+ G(x,t,D w, D^2 w)
 = 0
 & \mbox{in  $\R^N\times (0,T),$}\\[3mm]
 w (x,0)=\psi (x)& \mbox{in  $\R^N,$}\end{array}\right.\
\end{equation}
where $\Sigma (x,t) \in {\cal{M}}_N (\R)$ and $G$ is given by
(\ref{expression-G-intr}) and
\begin{equation}\label{hjb-related-bis}
\left\{\begin{array}{ll}
\displaystyle\frac{\partial  w}{\partial t}
+ h(x) |D w |^2 = 0 & \mbox{in  $\R^N\times (0,T),$}\\[3mm]
 w (x,0)=\psi (x)& \mbox{in  $\R^N,$}\end{array}\right.\
\end{equation}
where $h : \R^N \to \R.$
Our aim is to investigate comparison (and existence) results for
(\ref{hjb-related}) and  (\ref{hjb-related-bis})
under assumptions which include
Hamiltonians like (\ref{hamiltonien-related-1})
in Remark \ref{extension}(iii) (see Remark \ref{rmq-eq-rel}
for further comments).
More precisely,
we introduce two new assumptions:
\par\medskip
\noindent \textbf{(A4)}
$\Sigma\in C(\R^N\times [0,T]; {\cal{S}}_N^+ (\R))$ and,
for all $x\in\R^N, t\in [0,T],$
\begin{eqnarray*}
0< \Sigma(x,t) \ \ \ {\rm and} \ \ \ | \Sigma(x,t) | \leq \bar{C}.
\end{eqnarray*}
\noindent \textbf{(A5)}
$h\in C(\R^N ; \R),$
$h\in W^{2,\infty} (\tilde{\Gamma})$
where $\tilde{\Gamma}$ is an open neighborhood of
$\Gamma := \{ x\in \R^N : h(x)=0\}$
and, for all $x\in \Gamma,$
\begin{eqnarray*}
Dh(x) =0.
\end{eqnarray*}

\begin{Theorem} \label{thm-related}
Assume \textbf{(A2)}, \textbf{(A3)} and \textbf{(A4)}
(respectively \textbf{(A3)} and \textbf{(A5)}).
Let $U \in USC(\R^N\times [0,T])$ be a viscosity subsolution of
(\ref{hjb-related}) (respectively (\ref{hjb-related-bis}))
and $V \in LSC(\R^N\times [0,T])$ be a viscosity
supersolution of (\ref{hjb-related}) (respectively (\ref{hjb-related-bis}))
satisfying the
quadratic growth condition
(\ref{croiss-quad}).
Then $U\leq V$ in $\R^N\times [0,T].$
\end{Theorem}

The question of existence faces same problems as in
Section \ref{sec-1-comp}. We have existence and uniqueness
of a continuous viscosity solution for (\ref{hjb-related})
and (\ref{hjb-related-bis}) in the class of functions
with quadratic growth at least for short time
(as in Corollary \ref{thm-existence}). But solutions
can blow up in finite time.

Before giving the proof of the
theorem, we give some comments on the equations and the assumptions.

\begin{Remark} \rm \label{rmq-eq-rel} 
(i)
In the same way as in Remark \ref{extension}, above results hold
replacing $\Sigma$ by $-\Sigma$ in (\ref{hjb-related}) and $h$ by $-h$
in (\ref{hjb-related-bis}) and when dealing with terminal data
in both equations.  

(ii) Coming back to Remark \ref{extension}(iii), we note
$(a_1 +a_2 )(a_2-a_1)^T$ is not necessarily symmetric in
(\ref{hamiltonien-related-1}) whereas we assume $\Sigma$ to be symmetric
in \textbf{(A4)} (and $h$ is real-valued therefore symmetric
in \textbf{(A5)}). This is not a restriction
of generality since comparison for (\ref{hjb-related}) with
$\Sigma$ symmetric implies obviously comparison for
(\ref{hjb-related}) for any matrix $\Sigma,$
using that, for any $\Sigma\in {\cal{M}}_N (\R),$
$(\Sigma+\Sigma^T)/2 \in {\cal{S}}_N (\R).$ 

(iii) Coming back to Remark \ref{extension}(iv) again,
we see that \textbf{(A4)} corresponds to the case where the convex
Hamiltonian is predominant with respect to the concave one in
(\ref{hjbintr}) with (\ref{hamiltonien-related-1}).
Assumption \textbf{(A5)}
corresponds to (\ref{hamiltonien-related-1}) when $a_1$
and $a_2$ are real-valued but $h$ is allowed to change its sign.
For example, we have comparison for
\begin{eqnarray*}
 w_t + \phi(x)^3 |D w|^2 = 0,
\end{eqnarray*}
where $\phi \in C^2 (\R^N ;\R)$ is any bounded function.
Let us compare (\ref{eikonal-quad})
with first-order Hamilton-Jacobi equations whose Hamilonians are Lipschitz continuous both in the state
and gradient variables, or, to make simple, with the Eikonal equation
\begin{eqnarray*}
\displaystyle\frac{\partial  w}{\partial t}+a(x)|D w|=0 \ \ \ {\rm in} \ \R^N\times(0,T),
\end{eqnarray*}
where $a$ is Lipschitz continuous function. We know (\cite{cl83}, \cite{ley01}) we have existence and uniqueness
of a continuous viscosity solution for any continuous initial data without any restriction on the growth.
In this case, the sign of $a$ does not play any role whereas it seems to be the case for the
sign of $h$ in (\ref{eikonal-quad}). We would like to know if comparison for (\ref{eikonal-quad})
is true under weaker assumptions than \textbf{(A5)} and in dimension $N>1$ 
(i.e. when $\Sigma$ is neither positive nor negative definite in (\ref{hjb-related})).

(iv) There are some links between Theorem \ref{thm-unicite} and Theorem \ref{thm-related}.
Nevertheless we point out that, if Equation (\ref{hjbintr}) under assumptions of Section \ref{sec-1-comp}
is naturally associated with a control problem, this is not necessarily the case under the
assumptions of the current section. 
To be more precise, let us compare  (\ref{hjbintr}) with $H$ given by
(\ref{examplers}) and (\ref{hjb-related}) under \textbf{(A4)}. The
matrix $a$ can be singular in (\ref{examplers}) whereas we impose
the nondegeneracy condition $\Sigma>0$ in (\ref{hjb-related}). The
counterpart is that the regularity assumption with respect to $x$
on $\Sigma$ is weaker ($\Sigma$ is supposed to be merely continuous) than
the locally lipschitz regularity we assume for $a.$ Therefore
(\ref{hjb-related}) does not enter in the framework of Section
\ref{sec-1-comp}
in general.
A natural consequence is that proofs differ: in both  proofs of
Theorem \ref{thm-unicite} and \ref{thm-related}, the main argument
is a kind of linearization procedure, but while in the proof of
Theorem \ref{thm-unicite} we use essentially the convexity (or
concavity) of the operator corresponding to the unbounded control
set, here we use the locally Lipschitz continuity of the
Hamiltonian with respect to the gradient
uniformly in the state variable. 

(v) Note that we are able to deal with a second order term in
(\ref{hjb-related}) but not in (\ref{hjb-related-bis}).
\end{Remark}

{\it Proof of Theorem \ref{thm-related}.} We divide the
proof in two parts, corresponding respectively to the case of equation
\rec{hjb-related} and \rec{hjb-related-bis}. The main argument
in both proofs is a kind of linearization procedure as the one of
Lemma \ref{linearisation}. The rest of the proof is close to
the one of Theorem \ref{thm-unicite}. 

{\it Part 1.} We assume \textbf{(A2)}, \textbf{(A3)} and \textbf{(A4)}. \\
The estimates of $G$ are exactly the same than in
Lemma \ref{linearisation} so, for sake of simplicity, we
choose to take $G\equiv 0.$ The linearization procedure is
the subject of the following lemma whose proof is postponed
at the end of the section.
\begin{Lemma} \label{linearisation2}
Let $0<\lambda <1$ and set $\Psi= \lambda U-V.$ Then $\Psi$ is an $USC$
viscosity subsolution~of
\begin{equation}\label{equa-linearisee-2}
\left\{\begin{array}{ll}
\displaystyle\frac{\partial  w}{\partial t}
-\frac{2 \bar{C}}{1-\lambda}|D w|^2 \leq 0
& \mbox{in  $\R^N\times (0,T),$}\\[3mm]
 w (x,0)\leq (\lambda -1)\psi (x)& \mbox{in  $\R^N.$}
\end{array}\right.\
\end{equation}
\end{Lemma}
If $\psi \geq 0,$ then it follows $\Psi(x,0)\leq 0$ in $\R^N$ and using the
above lemma with the same supersolution as in the proof of
Theorem \ref{thm-unicite},
we obtain $\Psi\leq 0$ in $\R^N\times [0,T].$
Note that the boundedness of $\Sigma$ (see \textbf{(A4)}) is crucial
to build the supersolution. In the particular case when $G=0,$ we can also
choose a simpler supersolution such as for instance
\begin{eqnarray} \label{geant-sympa}
(x,t)\mapsto K \frac{[(|x|-R)^+]^2}{1-Lt}+\eta t, \ \ \ K,L,R,\eta >0
\end{eqnarray}
Letting $\lambda$ go  to 1, we conclude that  $U\leq V$ in $\R^N\times [0,T].$

It remains to prove that we can assume $\psi\geq 0$ without loss
of generality. To this end we use an argument similar to the one
of Step 2 in the proof of  Theorem \ref{thm-unicite}: let
$\bar{U}= U+M(1+|x|^2){\rm e}^{\rho t},$ $\bar{V}=
V+M(1+|x|^2){\rm e}^{\rho t},$ for some positive constants $M$ and
$\rho$ and set $\bar{\Psi}=\lambda \bar{U}-\bar{V}.$ We can easily
prove that, for $\rho > 16M^2\bar{C} {\rm e},$ $\bar{\Psi}$ is a
subsolution of (\ref{equa-linearisee-2}) (with a larger
constant $4 \bar{C}$ instead of $2\bar{C}$)
in $\R^N\times (0,\tau]$ for $\tau = 1/\rho.$
Moreover $\bar{\Psi} ({x},0)\leq (\lambda -1) (\psi(x) +
M(1+|x|^2){\rm e}^{\rho t})\leq 0$ for $M$ sufficiently large
(since $\psi$ has quadratic growth). Letting $\lambda$
go to 1, it follows that $U\leq V$ in $\R^N\times [0,\tau]$ and we conclude
by a step-by-step argument.

{\it Part 2.} We assume \textbf{(A3)} and \textbf{(A5)}. \\
Set $ \Omega^+ =\{ x\in\R^N : h(x)>0 \}$ and
$\Omega^- =\{ x\in\R^N : h(x)<0 \}$ which are open subsets of $\R^N.$
Define $\bar{\Psi}=\lambda \bar{U}-\bar{V}$ for $0<\lambda <1,$
$\bar{U}= U+M(1+|x|^2){\rm e}^{\rho t}$ and $\bar{V}=
V+M(1+|x|^2){\rm e}^{\rho t}$ with $M$ larger than the constants
which appear in  \textbf{(A3)} and in the growth condition
(\ref{croiss-quad}). Arguing as in Lemma \ref{linearisation2}
and noticing that all arguments are local, we obtain that
$\bar{\Psi}$ is a subsolution in $\Omega^+$ of the same kind of
equation as (\ref{equa-linearisee-2}), namely
\begin{equation}\label{equa-linearisee-333}
\left\{\begin{array}{ll}
\displaystyle\frac{\partial  w}{\partial t}
-\frac{C}{1-\lambda}|D w|^2 \leq 0
& \mbox{in  $ \Omega^+\times (0,\tau),$}\\[2mm]
 w (x,0)\leq 0& \mbox{in  $\Omega^+,$}
\end{array}\right.\
\end{equation}
for some constant $C=16\, ||h||_{\infty},$
$\rho > 4 MC{\rm e}$ and $\tau = 1/\rho >0$ which are independent
of $\lambda.$ Note that $\bar{\Psi}(\cdot, 0)\leq 0$ in
$\R^N$ because of the choice of $M.$

The sign of $\bar{\Psi}$ on $\Gamma$ is the subject of the following
lemma the proof of which uses \textbf{(A5)} and is postponed.
\begin{Lemma}\label{comp-sur-sigma}
For all $x\in\Gamma, \ t\in [0,T],$
we have $U(x,t)\leq \psi(x)\leq V(x,t).$
\end{Lemma}
Noticing that $\{ \Gamma , \Omega^+, \Omega^-\}$ is a partition of
$\R^N,$ we set
\begin{eqnarray*}
\hat{\Psi}:=
\left\{
\begin{array}{cc}
{\rm sup}\{\bar{\Psi},0\} & {\rm in} \  (\Gamma\cup \Omega^+)\times[0,T], \\
0                & {\rm in} \ \Omega^-\times[0,T]
\end{array}
\right.
\end{eqnarray*}
From the lemma, $\bar{\Psi}\leq 0$ in $\Gamma\times [0,T].$ Therefore
the function $\hat{\Psi}$ is continuous
in $\R^N\times (0,T).$ Moreover we claim that $\hat{\Psi}$ is a subsolution of
(\ref{equa-linearisee-333}) in $\R^N\times (0,\tau).$ This claim is clearly
true in $\Omega^-\times (0,T)$ (since $0$ is clearly a subsolution)
and in $\Omega^+\times (0,\tau)$ (since $\hat{\Psi}$ is
the supremum of two subsolutions). It remains to prove the result
on $\Gamma.$ Let $\varphi\in C^1 (\R^N\times (0,\tau))$ such that
$\hat{\Psi}-\varphi$ achieves a local maximum at a point $(\bar{x},\bar{t})
\in \Gamma\times (0,\tau).$ Let $(x,t)$ be in a neighborhood of
$(\bar{x},\bar{t}).$ If $x\in \Gamma\cup \Omega^+,$ then
$0\leq \hat{\Psi}(x,t)$ and if $x\in \Omega^-,$ then
$0= \hat{\Psi}(x,t).$ In any case,
$(0-\varphi)(x,t)\leq (\hat{\Psi}-\varphi)(x,t)
\leq (\hat{\Psi}-\varphi)(\bar{x},\bar{t}).$
But $\hat{\Psi}(\bar{x},\bar{t})\leq 0$ since $(\bar{x},\bar{t})
\in \Gamma\times (0,\tau).$
Therefore $(\bar{x},\bar{t})$ is a local maximum of $0-\varphi$
which ends the proof of the claim.
The initial condition $\hat{\Psi}\leq 0$ in $\R^N\times \{0\}$
is trivially satisfied. From the comparison principle proved
in Part 1 (for example using (\ref{geant-sympa}) as a supersolution),
we obtain $\hat{\Psi}\leq 0$ in $\R^N\times [0,\tau].$
Since $\tau$ does not depend on $\lambda,$ we can send $\lambda$ to $1.$
We obtain $U\leq V$ in $(\Gamma\cup \Omega^+) \times [0,\tau].$
We conclude in $(\Gamma\cup \Omega^+) \times [0,T]$ by a step-by-step
procedure.

Repeating the same kind of arguments replacing $\Omega^+$ by
$\Omega^-$ and $\bar{\Psi}$ by $U-M(1+|x|^2){\rm e}^{\rho t}
- \mu (V-M(1+|x|^2){\rm e}^{\rho t}),$ $0<\mu <1,$
we obtain that $U-V \leq 0$ in $(\Gamma\cup
\Omega^-)\times [0,T]$ which ends the proof.
~~$\Box$ 

{\it Proof of Lemma \ref{linearisation2}.}
We proceed as in the proof of Lemma \ref{linearisation}.
We can assume that $G\equiv 0$ since the computations with
$G$ are exactly the same than in Lemma \ref{linearisation}.

Note first that $\lambda U$ is an $USC$ subsolution of
\begin{equation}
\left\{\begin{array}{ll}
\displaystyle\frac{\partial  w}{\partial t}
+\frac{1}{\lambda} \langle \Sigma(x,t)D w ,D w \rangle \leq 0
& \mbox{in  $\R^N\times (0,T),$}\\[2mm]
 w (x,0)\leq \lambda \psi (x)& \mbox{in  $\R^N.$}
\end{array}\right.\
\end{equation}
Let $\varphi \in C(\R^N \times [0,T])$ and suppose that
$\Psi-\varphi$ reaches a strict local maximum at
$(\bar{x},\bar{t})\in \R^N\times (0,T]$ in some compact
subset $K\subset \R^N\times (0,T].$ We have
$$
\mathop{\rm max}_{(x,t), (y,t)\in K,}
\{ \lambda U(x,t)-{V} (y,t)-\varphi (\frac{x+y}{2},t)
-\frac{|x-y|^2}{\varepsilon^2} \}
\mathop{\longrightarrow}_{\varepsilon \downarrow 0}
\Psi(\bar{x},\bar{t})-\varphi(\bar{x},\bar{t}).
$$
The above maximum is achieved at some point
$(x_\varepsilon, y_\varepsilon, t_\varepsilon).$ Writing
the viscosity inequalities and subtracting them, we obtain
\begin{eqnarray*}
\varphi_t (\frac{x_\varepsilon +y_\varepsilon}{2},t_\varepsilon)
+ \frac{1}{\lambda}
\left\langle
\Sigma(x_\varepsilon, t_\varepsilon)
\left(\frac{1}{2}
D\varphi (\frac{x_\varepsilon +y_\varepsilon}{2},t_\varepsilon)
+ 2\frac{x_\varepsilon -y_\varepsilon}{\varepsilon^2}\right),
\frac{1}{2}
D\varphi (\frac{x_\varepsilon +y_\varepsilon}{2},t_\varepsilon)
+ 2\frac{x_\varepsilon -y_\varepsilon}{\varepsilon^2}
\right\rangle
\\ 
-\left\langle
\Sigma(y_\varepsilon, t_\varepsilon)
\left(-\frac{1}{2}
D\varphi (\frac{x_\varepsilon +y_\varepsilon}{2},t_\varepsilon)
+ 2\frac{x_\varepsilon -y_\varepsilon}{\varepsilon^2}\right),
-\frac{1}{2}
D\varphi (\frac{x_\varepsilon +y_\varepsilon}{2},t_\varepsilon)
+ 2\frac{x_\varepsilon -y_\varepsilon}{\varepsilon^2}
\right\rangle
\leq 0.
\end{eqnarray*}
In the sequel, we omit to write the dependence in $t_\varepsilon$ and
the point $((x_\varepsilon +y_\varepsilon)/2, t_\varepsilon)$
in the derivatives of $\varphi.$ Set
$p_\varepsilon = 2(x_\varepsilon -y_\varepsilon)/\varepsilon^2,$
$p_x = D\varphi /2 + p_\varepsilon$ and $p_y = -D\varphi /2 + p_\varepsilon.$
We have
\begin{eqnarray*}
0 &\geq &
\varphi_t + \frac{1}{\lambda} \langle \Sigma(x_\varepsilon)p_x , p_x \rangle
          - \langle \Sigma(y_\varepsilon) p_y, p_y \rangle \\
&\geq & \varphi_t +
\left( \frac{1}{\lambda} -1 \right) \langle \Sigma(x_\varepsilon)p_x , p_x \rangle
+  \langle (\Sigma(x_\varepsilon)-\Sigma(y_\varepsilon)) p_x, p_x \rangle
+\langle \Sigma(y_\varepsilon) p_x , p_x \rangle
-  \langle \Sigma(y_\varepsilon) p_y , p_y \rangle \\
&\geq & \varphi_t +
\left( \frac{1}{\lambda} -1 \right)
\langle \Sigma(x_\varepsilon)p_x , p_x \rangle
- m_K (|x_\varepsilon -y_\varepsilon|)|p_x|^2
+ \langle \Sigma(y_\varepsilon) p_x , p_x \rangle
-  \langle \Sigma(y_\varepsilon) p_y , p_y \rangle
\end{eqnarray*}
where $m_K$ is a modulus of continuity for $\Sigma$ in the compact
subset $K.$ Since $\Sigma(y_\varepsilon)$ is a symmetric matrix, we
have
\begin{eqnarray*}
\langle \Sigma(y_\varepsilon) p_x , p_x \rangle
-  \langle \Sigma(y_\varepsilon) p_y , p_y \rangle 
& = &
\langle \Sigma(y_\varepsilon ) (p_x+p_y) , p_x-p_y\rangle
=  2\, \langle \Sigma(y_\varepsilon) p_\varepsilon , D\varphi \rangle \\
& = & 
-\langle \Sigma(y_\varepsilon)D\varphi , D\varphi \rangle
+ 2\, \langle \Sigma(y_\varepsilon)p_x ,D\varphi \rangle.
\end{eqnarray*}
For any $\alpha >0,$ denoting by $\sqrt{\Sigma(y_\varepsilon)}$ the positive
symmetric squareroot of the positive symmetric matrix $\Sigma(y_\varepsilon),$
we get
\begin{eqnarray*}
 2\langle \Sigma(y_\varepsilon)p_x ,D\varphi \rangle
& = &
2\langle\sqrt{\Sigma(y_\varepsilon)}p_x ,\sqrt{\Sigma(y_\varepsilon)}D\varphi \rangle
\leq 
\alpha \langle \Sigma(y_\varepsilon)D\varphi ,D\varphi \rangle
+ \frac{1}{\alpha} \langle \Sigma(y_\varepsilon) p_x , p_x\rangle \\
&\leq & \alpha \langle \Sigma(y_\varepsilon)D\varphi ,D\varphi \rangle
+ \frac{1}{\alpha} \langle \Sigma(x_\varepsilon) p_x , p_x\rangle
+\frac{1}{\alpha}\langle m_K (|x_\varepsilon -
y_\varepsilon|)p_x,p_x\rangle.
\end{eqnarray*}
It follows
\begin{eqnarray*}
0 \geq \varphi_t + \left\langle \left[ \left( \frac{1}{\lambda} -1
- \frac{1}{\alpha}\right) \Sigma(x_\varepsilon) - \left(
1+\frac{1}{\alpha} \right) m_K (|x_\varepsilon - y_\varepsilon|)
Id \right] p_x , p_x \right\rangle -  (1+\alpha) \langle
\Sigma(y_\varepsilon)D\varphi ,D\varphi \rangle.
\end{eqnarray*}
Take $\alpha = \displaystyle\frac{2\lambda }{(1-\lambda)} >0$ in order to have
$\displaystyle\frac{1}{\lambda} -1 -\displaystyle\frac{1}{\alpha} =\displaystyle\frac{1}{2}
(\displaystyle\frac{1}{\lambda} -1) >0.$
We recall that $\Sigma>0$ by \textbf{(A4)},
hence for $\varepsilon$ small enough
the above scalar product is nonnegative. Thus
\begin{eqnarray*}
0\geq \varphi_t - (1+\alpha)  \langle
\Sigma(y_\varepsilon, t_\varepsilon)D\varphi , D\varphi\rangle
= \varphi_t - \frac{1+\lambda}{1-\lambda}  \langle
\Sigma(y_\varepsilon, t_\varepsilon)D\varphi , D\varphi\rangle.
\end{eqnarray*}
Letting $\varepsilon$ go to 0, we get
\begin{eqnarray*}
0\geq
\varphi_t(\bar{x},\bar{t}) - \frac{2}{1-\lambda}
\langle \Sigma(\bar{x},\bar{t}) D\varphi (\bar{x},\bar{t}),
D\varphi (\bar{x},\bar{t})\rangle
\end{eqnarray*}
which proves the result. ~~$\Box$

{\it Proof of Lemma \ref{comp-sur-sigma}.}
We make the proof for $U,$ the second inequality being similar.
Let $x_0\in \Gamma$ and consider, for $\eta >0,$
$$
\mathop{\rm sup}_{t\in [0,T]} \{ U(x_0,t) -\eta t\}.
$$
This supremum is achieved at a point $t_0\in [0,T]$ and we can assume,
up to subtract $|t-t_0|^2$ that it is a strict local maximum.
Consider, for $\varepsilon >0,$
$$
\mathop{\rm sup}_{\bar{B}(x_0,1)\times [0,T]}
\{ U(x,t) -\frac{|x-x_0|^2}{\varepsilon^2}-\eta t\}.
$$
This supremum is achieved at a point
$(x_\varepsilon,t_\varepsilon)$ and it is easy to see that
$(x_\varepsilon,t_\varepsilon) \to (x_0,t_0)$ when
$\varepsilon \to 0.$

Suppose that $t_\varepsilon >0.$ It allows us to
write the viscosity inequality for the subsolution $U$
at $(x_\varepsilon,t_\varepsilon)$ and  we get
\begin{eqnarray} \label{ineqvisco111}
\eta + h(x_\varepsilon)\left| 2\frac{x_\varepsilon-x_0}
{\varepsilon^2}\right|^2\leq 0.
\end{eqnarray}
Since $h\in W^{2,\infty}$ in a neighborhood of $\Gamma,$
we can write a Taylor expansion of $h$
at $x_0$ for $\varepsilon$ small enough
$$
h(x_\varepsilon) = h(x_0) + \langle Dh(x_0), x_\varepsilon -x_0\rangle
+\frac{1}{2} \langle D^2 h(x_0) (x_\varepsilon -x_0),
x_\varepsilon -x_0 \rangle + o(|x_\varepsilon -x_0|^2).
$$
From (\ref{ineqvisco111}) and \textbf{(A5)}, it follows
\begin{eqnarray*}
\eta -(2C+ m(|x_\varepsilon -x_0|^2))
\left(\frac{|x_\varepsilon -x_0|^2}{\varepsilon^2}\right)^2 \leq 0,
\end{eqnarray*}
where $m$ is a modulus of continuity.
Since ${|x_\varepsilon -x_0|^2}/{\varepsilon^2} \to 0$
as $\varepsilon\to 0,$ we obtain a contradiction
for small $\varepsilon.$

Therefore $t_\varepsilon =0$ for $\varepsilon$ small enough.
It follows that, for all $(x,t)\in \bar{B}(x_0,1)\times [0,T],$
we have
$$
U(x,t)-\frac{|x-x_0|^2}{\varepsilon^2}-\eta t \leq
U(x_\varepsilon, 0)-\frac{|x_\varepsilon-x_0|^2}{\varepsilon^2}
\leq \psi (x_\varepsilon).
$$
Setting $x=x_0$ and sending $\varepsilon$ and then $\eta$ to 0, we
obtain the conclusion. ~~$\Box$


\end{document}